\documentclass[journal]{IEEEtran}
\usepackage{cite}

\ifCLASSINFOpdf
\else
\fi
\usepackage{amsmath}
\usepackage{amsfonts}
\usepackage{amssymb}

\usepackage{amsmath}
\usepackage{booktabs}
\usepackage{tikz}

\usetikzlibrary{backgrounds,calc, arrows, intersections}
\tikzset{
  sin v source/.style={
    circle,
    draw,
    append after command={
      \pgfextra{
      \draw
      ($(\tikzlastnode.center)!0.5!(\tikzlastnode.west)$)
       arc[start angle=180,end angle=0,radius=0.425ex]
      (\tikzlastnode.center)
       arc[start angle=180,end angle=360,radius=0.425ex]
      ($(\tikzlastnode.center)!0.5!(\tikzlastnode.east)$)
    ;
    }
  },
  scale=1.5,
 }
}

\usetikzlibrary{plotmarks}
\usetikzlibrary{calc}
\usetikzlibrary{arrows}
\usetikzlibrary{arrows,intersections}
\usetikzlibrary{arrows.meta}
\usetikzlibrary{shapes,decorations}
\usetikzlibrary{positioning}
\usepackage{pgfplots}
\tikzset{sin v source/.style={
  circle,
  draw,
  append after command={
    \pgfextra{
    \draw
      ($(\tikzlastnode.center)!0.5!(\tikzlastnode.west)$)
       arc[start angle=180,end angle=0,radius=0.425ex]
      (\tikzlastnode.center)
       arc[start angle=180,end angle=360,radius=0.425ex]
      ($(\tikzlastnode.center)!0.5!(\tikzlastnode.east)$)
    ;
    }
  },
  scale=1.5,
 }
}

\pgfplotsset{compat=1.14}
\usepgfplotslibrary{statistics}

\usepackage{circuitikz}
\usepackage{tikz}
\usepackage{placeins}

\newcommand*{\bl}{\textcolor{black}}

\newcommand*{\blue}{\textcolor{black}}

\newcommand*{\blueEPSR}{\textcolor{black}}

\newcommand*{\rev}{\textcolor{black}}

\newcommand*{\red}{\textcolor{black}}

\begin{document}

\title{Inexact Convex Relaxations for AC Optimal \\ Power Flow: Towards AC Feasibility}

\author{Andreas~Venzke, 
        Spyros~Chatzivasileiadis, 
        and~Daniel~K.~Molzahn \vspace{-0.3cm}
\thanks{A. Venzke and S. Chatzivasileiadis are with the Department of Electrical Engineering, Technical University of Denmark, 2800 Kgs. Lyngby, Denmark e-mail: \{andven, spchatz\}@elektro.dtu.dk.}
\thanks{D. K. Molzahn is with the School of Electrical and Computer Engineering, Georgia Institute of Technology, Atlanta, GA 30313, USA as well as the Energy Systems Division, Argonne National Laboratory, Lemont, IL 60439 USA, e-mail: molzahn@gatech.edu.}
\thanks{This work  is  supported  by  the  multiDC  project,  funded  by  Innovation  Fund Denmark, Grant Agreement No. 6154-00020B as well as by the U.S. Department of Energy, Office of Electricity Delivery and Energy Reliability under contract DE-AC-02-06CH11357.}
}

\maketitle
\begin{abstract}
Convex relaxations of AC optimal power flow (AC-OPF) problems have attracted significant interest as in several instances they provably yield the global optimum to the original non-convex problem. If, however, the relaxation is inexact, the obtained solution is not AC-feasible. The quality of the obtained solution is essential for several practical applications of AC-OPF, but detailed analyses are lacking in existing literature. This paper aims to cover this gap. We provide an in-depth investigation of the solution characteristics when convex relaxations are inexact, we assess the most promising AC feasibility recovery methods for large-scale systems, and we propose two new metrics that lead to a better understanding of the quality of the identified solutions. We perform a comprehensive assessment on 96 different test cases, ranging from 14 to 3120 buses, and we show the following: (i)~Despite an optimality gap of less than 1\%, several test cases still exhibit substantial distances to both AC feasibility and local optimality and the newly proposed metrics characterize these deviations. (ii)~Penalization methods fail to recover an AC-feasible solution in 15 out of 45 cases. (iii)~The computational benefits of warm-starting non-convex solvers have significant variation, but a computational speedup exists in over 75\% of the cases. 
\end{abstract}


\IEEEpeerreviewmaketitle
\section{Introduction} 
The AC optimal power flow (AC-OPF) problem is \rev{fundamental} for the efficient operation of power systems~\cite{Panciatici_PSCC}. \red{Formulations of AC-OPF have found practical use in tools that minimize system losses and optimize setpoints of reactive power sources (e.g., synchronous condensers). Moreover, \mbox{AC-OPF} is being increasingly considered for market clearing procedures.} 
The \mbox{AC-OPF} minimizes an objective function (e.g., generation cost) subject to the power system  operational constraints (e.g., limits on the transmission line flows and bus voltages). However, nonlinearities from the AC power flow equations result in the AC-OPF problem being non-convex and generally NP-hard~\cite{Hentenryck_2016_NP,bienstock2019strong}.
To address this issue, different convex relaxations of the AC-OPF problem have been proposed during the last decade, including second-order cone programming (SOCP)~\cite{jabr2006radial}, semidefinite programming (SDP)~\cite{bai2008semidefinite, Lavaei2012}, and quadratic convex (QC) relaxations~\cite{coffrin2016qc}. These relaxations have attained significant interest as in several test cases they provably yield the global optimum to the original non-convex problem~\cite{Lavaei2012,molzahn2013implementation}, i.e., they are exact, and are shown to be tractable for test cases with thousands of buses~\cite{molzahn2013implementation,coffrin2016qc}. Besides obtaining a global optimality certificate, solving a convex instead of a non-convex problem has major advantages in various applications of AC-OPF. For example, several decomposition techniques are only guaranteed to converge for convex problems~\cite{molzahn2017survey} and bi-level programs arising in AC-OPF under uncertainty are more tractable~\cite{lorca2017robust}. \bl{Additional examples include applications of convex relaxations for distributed AC-OPF in microgrids \cite{dall2013distributed}, and for reactive power control in distribution networks \cite{zheng2015fully}.}

If, however, the convex relaxations of AC-OPF problems are inexact, the obtained solutions are no longer AC-feasible, yielding only a lower bound on the objective value. 
This poses a barrier for practical applications. Understanding when convex relaxations fail to be exact and what are the most promising options to obtain an AC-feasible \mbox{(near-)globally optimal} solution becomes fundamental for enabling the use of these methods in practice. This paper aims to cover this gap. To this end, we provide an in-depth analysis of the solution characteristics when convex relaxations are inexact, we assess promising AC feasibility recovery methods in a wide range of cases, and we propose new metrics that lead to a better understanding of the quality of the identified solutions. 

\bl{A related strand in research investigates theoretical conditions which guarantee exactness of convex relaxations of AC-OPF problems. These are, however, mainly limited to radial networks (see, e.g., \cite{low2014convex} for a comprehensive review) or impose restrictive assumptions on network parameters \cite{mahboubi2018analysis}.}
\bl{On the other hand, there are several works \cite{kocuk2015inexactness,louca2014nondegeneracy} which define sufficient conditions for inexactness of convex relaxations. The work in \cite{kocuk2015inexactness} focuses on radial networks, and explores conditions for which the SDP relaxation is inexact. The work in \cite{louca2014nondegeneracy} provides conditions for arbitrary networks under which the solution to the SDP relaxation will have rank larger than one, i.e., be inexact. The motivation for our work is that the derived conditions on exactness and inexactness do not apply to a large set of meshed transmission grid test instances. Instead, using an comprehensive empirical analysis, we investigate the solution quality when relaxations are inexact, define two metrics to evaluate the distance to AC-feasibility and local optimality, and evaluate the most promising AC-feasibility recovery procedures.} 

First, we assess the quality of the identified solutions when convex relaxations are inexact. \rev{Previous research has shown that although convex relaxations are inexact for the majority of available test cases~\cite{PGLIB}, optimality gaps (i.e., the difference between the objective value of the relaxation and the objective value reported by a non-convex solver) of less than 1\% can be achieved in many instances. While much of the literature (e.g., \cite{coffrin2016qc, kocuk2016strong}) focuses on further reducing the optimality gap, the quality of the obtained solution is often neglected.} However, the quality of the obtained decision variables is essential for various applications of AC-OPF, e.g., in bi-level programs~\cite{lorca2017robust} or where an AC-feasible solution is a requirement. As we show in this paper, it is important to realize that even a zero optimality gap does not guarantee that the obtained solution is AC-feasible.

\rev{In this paper, we provide a comprehensive assessment of the inexact solutions to the QC and SDP relaxations with respect to both AC feasibility and local optimality. From our analysis, it becomes obvious that the optimality gap alone is an insufficient metric for assessing the quality of the obtained solution. To address this limitation, we propose two new metrics: i)~the cumulative normalized constraint violation, and ii)~the average normalized distance to local optimality. While existing related studies have focused only on radial distribution networks and the obtained voltage magnitudes for inexact SOCP relaxations~\cite{abdelouadoud2015optimal}, in our comprehensive assessment, we use a wide range of \rev{meshed transmission network} PGLib OPF test cases from~\cite{PGLIB} and consider all AC-OPF state variables.}

\rev{Second, we rigorously evaluate two of the most promising directions for recovering an AC-feasible solution on up to 96 different test cases, ranging from 14 up to 3120 buses. The first focuses on modifying the objective functions of convex relaxations with penalization terms to guide them towards an AC-feasible solution~\cite{Madani2015, Madani2016,Wei2017,liu2019rank}. \bl{We show how the newly proposed metrics can assess and improve the performance of penalization methods.}}

\rev{This paper focuses on penalty terms based on reactive power~\cite{Madani2015} and apparent branch flow losses~\cite{Madani2016}, as they have been shown to be tractable for larger systems and to result in near-globally optimal solutions for certain test cases.} 

The second direction to recover an AC-feasible solution uses the result of the inexact convex relaxation to warm-start a general non-linear solver. While prior work in~\cite{kocuk2016strong} has solely focused on warm-starting interior-point solvers with solutions to the inexact SOCP relaxation, this paper investigates warm-starting both interior-point and sequential quadratic programming solvers with inexact solutions to the QC and SDP relaxations. 

The main contributions of this paper are:
\begin{enumerate}
    \item This is the first work to provide an in-depth assessment of the quality of the solution obtained through convex relaxations, measuring the distance of the decision variables for both the QC and SDP relaxations to AC feasibility and local optimality in 96 test cases, ranging from 14 to 3120 buses. We propose two empirical metrics complementary to the optimality gap: i)~the cumulative normalized constraint violation, and ii)~the average normalized distance to local optimality. We show that despite an optimality gap of less than 1\%, several test cases still exhibit substantial distances to AC feasibility and local optimality, highlighting the added value of the two metrics. 
    
    \item We provide a rigorous analysis of three different penalization methods for the SDP relaxation on 45 PGLib OPF test cases \rev{up to 300 buses. We show that they fail to recover an AC-feasible solution for 35.6\% of test cases and can incur significant sub-optimality of up to 47.7\%.} \bl{We characterize the obtained solutions from penalized SDP relaxations using our proposed metrics. We also show that in cases where penalization methods fail, they often exhibit substantial distances to both AC-feasibility and local optimality. For failed test instances with small distances, we show how our proposed metrics inform a fine-tuning of penalty weights to obtain AC-feasible solutions.}
    
    \item \bl{We investigate warm-starting interior-point and sequential quadratic programming solvers with the solutions to the inexact convex relaxations compared to initializations with a flat start or the solutions to the DC optimal power flow}. Examining 96 test cases with up to 3120 buses, we show that benefits in terms of computational speed, solver reliability, and solution quality strongly depend on solver and test case, which corroborates the complexity of the AC feasibility recovery problem. The warm-started interior point solver IPOPT \cite{wachter2006implementation} achieves the best performance, gaining a computational speed-up in over 75\% of the cases.

\end{enumerate}

This paper is structured as follows: Section~\ref{sec:formulation} formulates the AC-OPF problem and the considered QC and SDP relaxations. Section~\ref{Metrics} proposes two metrics to assess the distances to AC feasibility and local optimality of inexact convex relaxations. Section~\ref{sec:recovery} reviews different methods for recovering AC-feasible or locally optimal solutions from inexact convex relaxations. Section~\ref{sec:results} provides extensive computational studies using the PGLib OPF benchmarks. Section~\ref{sec:conclusion} concludes.
\section{AC Optimal Power Flow and Relaxations}
\label{sec:formulation}

The AC-OPF problem has a variety of different mathematically equivalent formulations. For a detailed survey on AC-OPF and convex relaxations, the reader is referred to \cite{molzahn2018fnt}. Here, for brevity, we follow the AC-OPF formulation of \cite{coffrin2016qc} to facilitate the derivation of the SDP and QC relaxations. 

A power grid consists of the set $\mathcal{N}$ of buses, a subset of those denoted by $\mathcal{G}$ have a generator. The buses are connected by a set $(i,j) \in \mathcal{L}$ of power lines from bus $i$ to $j$. The optimization variables are the complex bus voltages $V_k$ for each bus $k \in \mathcal{N}$ and the complex power dispatch of generator $S_{G_k}$ for each bus $k \in \mathcal{G}$. The objective function $f_{\text{cost}}$ minimizes the cost associated with active power dispatch:
\begin{align}
\min_{V,{S_G}} \, f_{\text{cost}}:= \sum_{k \in \mathcal{G}} c_{k_2} \Re \{S_{G_k}\}^2 + c_{k_1} \Re \{S_{G_k}\} + c_{k_0} \label{Obj}
\end{align}
The terms $c_{k_2}$, $c_{k_1}$ and $c_{k_0}$ denote quadratic, linear and constant cost terms associated with generator active power dispatch, respectively. The following constraints are enforced:
\begin{subequations}
\label{opf}
\begin{align}
  &  (V_k^{\text{min}})^2 \leq V_k V_k^* \leq (V_k^{\text{max}})^2 & \quad  \forall k \in \mathcal{N} \label{Vmax} \\
  &  S_{G_k}^{\text{min}} \leq S_{G_k} \leq S_{G_k}^{\text{max}} & \quad  \forall k \in \mathcal{G} \label{Smax} \\
  &  |S_{ij}| \leq S_{ij}^{\text{max}} & \quad \forall (i,j) \in \mathcal{L} \label{Sijmax} \\
  &  S_{G_k} - S_{D_k} = \sum_{(k,j) \in \mathcal{L}} S_{kj} & \quad \forall k \in \mathcal{N} \label{Sbal} 
     \end{align}
  \begin{align}
  & S_{ij} = Y_{ij}^* V_i V_i^* - Y_{ij}^*V_i V_j^* & \quad \forall (i,j) \in \mathcal{L} \label{Sij} \\
  & -\theta^{\text{max}}_{ij} \leq \angle (V_i V_j^*) \leq \theta^{\text{max}}_{ij} & \quad \forall (i,j) \in \mathcal{L} \label{angle} 
\end{align}
\end{subequations}
The bus voltage magnitudes are constrained in \eqref{Vmax} by upper and lower limits $V_k^{\text{min}}$ and $V_k^{\text{max}}$. The superscript $*$ denotes the complex conjugate. Similarly, the generators' complex power outputs are limited in \eqref{Smax} by upper and lower bounds $S_{G_k}^{\text{min}}$ and $S_{G_k}^{\text{max}}$, where inequality constraints for complex variables are interpreted as bounds on the real and imaginary parts. The apparent branch flow $S_{ij}$ is upper bounded in \eqref{Sijmax} by $S_{ij}^{\text{max}}$. The nodal complex power balance \eqref{Sbal} including the load $S_D$ has to hold for each bus. The apparent branch flow $S_{ij}$ is defined in \eqref{Sij}. The term $Y$ denotes the admittance matrix of the power grid.  The branch flow is also limited in \eqref{angle} by an upper limit on the angle difference $\theta^{\text{max}}_{ij}$. As proposed in \cite{Lavaei2012, coffrin2016qc}, an additional auxiliary matrix variable $W$ is introduced, which denotes the product of the complex bus voltages:
\begin{align}
    W_{ij} = V_i V_j^* \label{Wnonconvex}
\end{align}
This facilitates the reformulation of \eqref{Vmax}, \eqref{Sij}, and \eqref{angle} as:
\begin{subequations}
\label{Wconstraints}
\begin{align}
      &  (V_k^{\text{min}})^2 \leq W_{kk} \leq (V_k^{\text{max}})^2 & \quad  \forall k \in \mathcal{N} \label{Wmax} \\
       & S_{ij} = Y_{ij}^* W_{ii} - Y_{ij}^* W_{ij} & \quad \forall (i,j) \in \mathcal{L} \label{WSij} \\
           & S_{ij} = Y_{ij}^* W_{jj}- Y_{ij}^* W_{ij}^*  & \quad \forall (i,j) \in \mathcal{L} \label{WSji} \\
             & \tan(-\theta^{\text{max}}_{ij}) \leq \tfrac{\Re \{W_{ij}\}}{\Im \{W_{ij}\}} \leq \tan(\theta^{\text{max}}_{ij}) & \quad \forall (i,j) \in \mathcal{L} \label{Wangle} 
\end{align}
\end{subequations}
The only source of non-convexity is the voltage product \eqref{Wnonconvex}.
\subsection{DC Optimal Power Flow (DC-OPF) Approximation}
The DC-OPF, which serves here as benchmark, is an approximation that is often used in, e.g., electricity markets and unit commitment problems. This approximation neglects voltage magnitudes, reactive power, and active power losses. The state variables are the active generation $P_G$ and the voltage angles $\theta$. Depending on whether a quadratic cost term is included, the optimization problem is either a linear program (LP) or a quadratic program (QP). For brevity the formulation is omitted, but can be found in, e.g., \cite{dvijotham_molzahn}. \rev{Since solver reliability is very important, most system operators currently use the DC power flow approximation, which, however, can exhibit substantial errors~\cite{dvijotham_molzahn}.} 
\subsection{Quadratic Convex (QC) Relaxation}
The QC relaxation proposed in \cite{coffrin2016qc} uses convex envelopes of the polar representation of the AC-OPF problem to relax the dependencies among voltage variables. The non-convex constraint \eqref{Wnonconvex} is removed; variables for voltages, $v_i \angle \theta_i \, \,\forall i \in \mathcal{N}$, and squared current flows, $l_{ij} \, \,\forall (i,j) \in \mathcal{L}$, are added; and convex constraints and envelopes are introduced~\cite{coffrin2016qc}:
\begin{subequations}
\small
\label{QC}
\begin{align}
  & W_{kk} = \left\langle v_k^2 \right\rangle^T  &  \forall k \in \mathcal{N} \label{QC1}\\
  & \Re \{W_{ij}\} =   \left\langle \left\langle v_i v_j \right\rangle^M   \left\langle \cos (\theta_i - \theta_j) \right\rangle^C \right\rangle^M  &  \forall (i,j) \in \mathcal{L}\\
  & \Im \{W_{ij}\} =   \left\langle \left\langle v_i v_j \right\rangle^M   \left\langle \sin (\theta_i - \theta_j) \right\rangle^S \right\rangle^M &   \forall (i,j) \in \mathcal{L}\\
  & S_{ij} + S_{ji} = Z_{ij} l_{ij} &  \forall (i,j) \in \mathcal{L}\\
  & |S_{ij}|^2 \leq W_{ii} l_{ij} &  \forall (i,j) \in \mathcal{L} \label{QCend}
\end{align}
\end{subequations}
The superscripts $T,M,C,S$ denote convex envelopes for the square, bilinear product, cosine, and sine functions,  respectively; for details, see \cite{coffrin2016qc}. The term $Z_{ij}$ denotes the line impedance. The resulting optimization problem is an SOCP that minimizes \eqref{Obj} subject to \eqref{Smax} -- \eqref{Sbal}, \eqref{Wconstraints}, and \eqref{QC}. The QC relaxation dominates the SOCP relaxation of~\cite{jabr2006radial} in terms of tightness at a similar computational complexity, and is particularly effective for meshed transmission networks with tight angle constraints. We therefore omit the SOC relaxation. 
\subsection{Semidefinite (SDP) Relaxation}
In the semidefinite (SDP) relaxation proposed in \cite{bai2008semidefinite, Lavaei2012}, the non-convex constraint \eqref{Wnonconvex} is reformulated in matrix form:
\begin{subequations}
\begin{align}
    W \succeq 0 \label{SDP} \\
    \text{rank}(W) = 1 \label{rank} 
\end{align}
\end{subequations}
The non-convexity of the resulting formulation is encapsulated in the rank constraint \eqref{rank}, which is subsequently relaxed. The resulting optimization is an SDP that minimizes the objective function \eqref{Obj} subject to \eqref{Smax} -- \eqref{Sbal}, \eqref{Wconstraints}, and \eqref{SDP}. In terms of theoretical tightness, the QC neither dominates nor is dominated by the SDP relaxation \cite{coffrin2016qc}. 

\bl{Higher-order moment relaxations generalize the SDP relaxation of~\cite{Lavaei2012}, facilitating global solutions to a broader class of problems at the computational cost of larger SDP constraints~\cite{josz2014application, molzahn2014moment}. To reduce the size of these constraints, the work in \cite{molzahn2015sparsity} proposes a method to exploit sparsity by sequentially enforcing higher-order moment constraints only for parts of the network which corresponds to higher rank-components in matrix $W$. As this approach is not as computationally mature as the SDP and QC relaxations, and its use for various further application of the AC-OPF, e.g. distributed AC-OPF and bi-level programs, is not explored yet, we leave its detailed computational analysis for future work.}

\section{Distance Metrics for Inexact Solutions}
\label{Metrics}
\blue{The optimality gap is a widely used distance metric for assessing the quality of inexact solutions obtained from convex relaxations. The optimality gap between the solution to the convex relaxation and the best known feasible point for the non-convex AC-OPF problem is defined as follows:
\begin{align}
    (1 - \tfrac{f_{\text{cost}}^{\text{relax}}}{f_{\text{cost}}^{\text{local}}}) \times 100 \%
\end{align}
The term $f_{\text{cost}}^{\text{relax}}$ denotes the lower objective value bound from the relaxation and $f_{\text{cost}}^{\text{local}}$ is the objective value of the best known feasible point obtained from a local non-convex solver.  If the relaxation is inexact, the magnitude of the optimality gap does not necessarily indicate the decision variables' distances to feasibility or local optimality for the original non-convex AC-OPF problem; e.g., for cases with very flat objective functions, solutions with small optimality gaps could still exhibit substantial distances to both. Additionally, the closest solution that is AC-feasible might not coincide with the closest locally optimal solution. To assess both these distances for a wide range of test cases, we propose two new alternative metrics: i)~the cumulative normalized constraint violation and ii)~the average normalized distance to a local solution.}
\vspace{-0.2cm}
\subsection{Cumulative Normalized Constraint Violation}
\label{Metrics:constrviol}
To assess the distance to AC feasibility, we run an AC power flow (AC-PF) with set-points obtained from an inexact convex relaxation and evaluate the constraint violations. For each bus~$k\in\mathcal{N}$, there are four state variables: the voltage magnitude $|V_k|$, the voltage angle $\theta_k$, the active power injection $\Re\{S_{Gk}-S_{Dk}\}$, and the reactive power injection $\Im \{S_{Gk}-S_{Dk}\}$. 
If the solution to the convex relaxation is inexact, the resulting state variables do not fulfill the AC power flow equations and are hence AC-infeasible. To recover a solution that satisfies the AC power flow equations, an AC power flow (AC-PF) can be computed with the set-points obtained from the convex relaxation. 

Theoretically, it is possible to model \blueEPSR{generator-connected buses} as $PV$, $PQ$, or $V\theta$ buses. Please note that for all three options discussed in the following we assume that one of the generator is designated as slack bus, and is therefore modelled as $V\theta$ bus. We exclude the possibility to model all generators as $V\theta$ buses due to two reasons: First, the voltage angles are not readily available for the SDP relaxation, and, if estimated, would likely of lower quality than $P$, $Q$ and $V$ set-points directly obtained from the solution to the relaxation. Second, specifying voltage angles at multiple buses within the AC power flow has to our knowledge not been reported in literature, and we expect poor convergence performance. For modeling as $PQ$ buses, in our computational experiments, we also observed issues with the convergence of AC power flows. Accordingly, we choose to use the setpoints for active power injections and voltage magnitudes at \blueEPSR{generator-connected buses}  from the relaxation's solution, i.e., modeling generators as $PV$ buses. 

We propose the cumulative normalized constraint violation resulting from this power flow solution as a metric to quantify the distance to AC feasibility. This metric can be computed by taking the sum of constraint violations from the AC-PF solution, normalized by the respective upper and lower constraint bounds. \rev{For each variable $x:=\{P_G, Q_G, |V|, \theta_{ij}, S_{ij}\}$ and corresponding set $\mathcal{X} =\{\mathcal{G},\mathcal{G},\mathcal{N},\mathcal{L},\mathcal{L}\}$ we define the cumulative normalized constraint violations $x_{\text{viol}}$ as:}
\begin{align}
  \rev{  x_{\text{viol}} }&\rev{:= \sum_{k \in \mathcal{X}} \tfrac{\max(x_k^{\text{AC-PF}}-x_k^{\text{max}},x_k^{\text{min}}-x_k^{\text{AC-PF}},0)}{x_k^{\text{max}}-x_k^{\text{min}}}\times100\% } \label{eq:constrviol} 
\end{align}
\blueEPSR{The term $x_k^{\text{AC-PF}}$ refers to to the value of the variable from the AC power flow solution. The minimum and maximum constraints on the variable are represented by the terms $x_k^{\text{max}}$ and $x_k^{\text{min}}$, respectively.}
 \bl{We take the magnitude of the violation into account since larger violations indicate a larger distance to AC-feasibility. We normalize each violation with respect to the upper and lower constraint limits, which allows for a fair comparison between different type of violations. We choose not to normalize by the number of buses or constraints, as only a very low percentage of constraints is usually violated in our studies, as can be seen in Fig.~\ref{Mean_Max_Viol}.} We found that this metric carries more information than assessing the average or maximum constraint violations since only a small subset of the constraints are active in typical AC-OPF problems. This metric allows for a comparison of the distance to AC feasibility among relaxations for a given system. Note that it is not averaged by the number of buses, and as a result, inexact convex relaxations of AC-OPF problems for larger systems may exhibit larger values.

\blueEPSR{We note that the accuracy of any power system model inherently relies on appropriate values for parameters such as transmission line limits, impedances, etc. These parameters can vary with changes in the weather and operating conditions, leading to differences between the modeled and actual system behavior. While analyses of these differences are not in the scope of this paper, we refer the reader to other research which proposes models that incorporate the effects of temperature, weather, frequency variations, etc. on the model parameters~\cite{RAHMAN2019241, POUDEL2017266, 8351381, 8274733}.}

\subsection{Average Normalized Distance to a Local Solution} \label{Metrics:local}
To assess the distance to local optimality, we propose an additional metric defined as the averaged normalized distance of the variable values obtained with the inexact convex relaxation to a locally optimal solution obtained with a non-convex solver.\footnote{\rev{Note that no AC-OPF problems with multiple local solutions were identified in the numerical results for our large-scale test cases with any of the three different solvers or in~\cite{kardos2018complete}. This could indicate that no spurious locally optimal solutions exist, i.e., solutions that are locally optimal and not globally, and the identified solutions are in fact the globally optimal solutions.}} This metric can be computed by taking the absolute difference between the variable value from the relaxation's solution and the value from the local optimum, normalized by the difference in the upper and lower variable bound, and then averaging over all variables. As this metric is averaged by the number of variables, which is a function of the system size, it allows for fair comparison among relaxations for systems with different sizes. \rev{For each variable $x:=\{P_G, Q_G, |V|, \theta_{ij}, S_{ij}\}$ and corresponding set $\mathcal{X} =\{\mathcal{G},\mathcal{G},\mathcal{N},\mathcal{L},\mathcal{L}\}$, we define the average normalized distance $x_{\text{dist}}$ as:} 
\begin{align}
      \rev{ x_{\text{dist}}} &\rev{: = \tfrac{1}{|\mathcal{X}|} \sum_{k \in \mathcal{X}} \tfrac{| x_k^{\text{relax}} - x_k^{ \text{opt}} |}{x_k^{\text{max}}-x_k^{\text{min}}} \times 100\%} \label{eq:local}
\end{align}
\blueEPSR{The terms $x_k^{\text{relax}}$ and $x_k^{ \text{opt}}$ are the values of the variable defined by the solution to the convex relaxation and by a locally optimal solution to the original non-convex AC-OPF problem, respectively.}

\section{Recovery of AC-Feasible and Local Solutions}
\label{sec:recovery}

\subsection{Penalization Methods for SDP Relaxation}
\label{Theory_Pen}
For a variety of test cases, the SDP relaxation has a small optimality gap but the obtained solution $W$ does not fulfill the rank-1 condition \eqref{rank}; see, e.g., \cite{Madani2015}. To drive the solution towards a rank-1 point in order to recover an AC-feasible solution, several works~\cite{Madani2015,Madani2016} have proposed augmenting the objective function of the semidefinite relaxation \blueEPSR{$f_{\text{cost}}$} with penalty terms, denoted with $f_{\text{pen}}$. Note that the penalized formulations \emph{are not} relaxations of the original AC-OPF problem, but can still be useful for recovering feasible points. We focus on the following three penalty terms. First, the nuclear norm proposed by \cite{fazel2001rank} is a widely used penalty term for the general rank-minimization problem: 
 \begin{align}
     f_{\text{pen}} = f_{\text{cost}} + \epsilon_{\text{pen}} \text{Tr} \{ W \} \label{T_Pen}
 \end{align}
 The penalty weight is denoted with $\epsilon_{\text{pen}}$ and the term $\text{Tr} \{\,\cdot\, \}$ indicates the trace of the matrix $W$, i.e., the sum of the diagonal elements. Second, specific for the AC-OPF, the work in \cite{Madani2015} proposes penalization of the reactive generator outputs:
 \begin{align}
     f_{\text{pen}} = f_{\text{cost}} + \epsilon_{\text{pen}}\sum_{k \in \mathcal{G}} \Im \{S_{G_k}\} \label{Q_Pen}
 \end{align}
 \blueEPSR{The term $\Im \{S_{G_k}\}$ denotes the reactive power output of the generator connected to bus $k$.} Finally, the work in \cite{Madani2016} suggests adding an apparent branch flow loss penalty to the objective function: 
  \begin{align}
      f_{\text{pen}} = f_{\text{cost}} + \epsilon_{\text{pen}} \sum_{(i,j) \in \mathcal{L}} |S_{ij} - S_{ji}| \label{S_Pen}
        \vspace{-0.2cm}
  \end{align}
\blueEPSR{The term $|S_{ij} - S_{ji}|$ represents the absolute apparent power loss on the line from bus $i$ to bus $j$.} Both the works \cite{Madani2015, Madani2016} show that certain choices of $\epsilon_{\text{pen}}$ result in successful recovery of an AC-feasible and near-globally optimal operating points for selected test cases. In Section~\ref{Pen_results}, we will present a counterexample and provide a detailed empirical analysis of the different penalty terms for a wide range of test cases. Note that to the knowledge of the authors, penalization methods have not been applied to address inexact solutions resulting from the QC relaxation.

\bl{The two proposed metrics in Section~\ref{Metrics:constrviol} and Section~\ref{Metrics:local} complement the evaluation of penalization methods as follows. First, by evaluating the cumulative normalized constraint violation for different penalty weights $\epsilon_{\text{pen}}$, we can quantify if a solution obtained from the penalized SDP relaxation is AC-feasible as well as how different magnitudes of penalty weights compare in terms of distance to AC-feasibility, e.g., we can assess whether increasing the penalty moves the obtained solution closer to AC-feasibility. Second, by evaluating the average normalized distance to local optimality, we can assess whether the penalization terms drive the solution towards the solution provided by a non-convex solver, or otherwise either towards regions with sub-optimal costs. Note that in theory the penalization could also drive the solution towards a spurious locally or globally optimal solution with an lower objective value. For the large-scale test cases in our simulation study, we did not observe multiple locally optimal solutions, and the penalization did not drive the solution towards an AC-feasible solution with lower objective value. Third, if both metrics remain substantial for a wide range of penalty weights $\epsilon_{\text{pen}}$, then this indicates that the penalty term under study might not be able to recover an AC-feasible solution. Refer to Section~\ref{Pen_results} for a detailed computational analysis.}

\vspace{-0.2cm}
\subsection{Warm-Starting Non-Convex Local Solvers} \label{Warm_start_theory}
When penalization methods are not successful at recovering a rank-1 solution, non-convex solvers can be warm started with the solution of convex relaxations in order to recover an AC-feasible and locally optimal solution. Compared to a flat start of $V_k = 1\angle 0$, $\forall k\in\mathcal{N}$, which is a common initialization for non-convex solvers, warm-starting could lead to i)~reduced computational time and ii)~improved solution quality. For these purposes, we utilize two types of non-convex solvers:
\subsubsection{Sequential Quadratic Programming (SQP)} To compute a search direction, this method iteratively solves second-order Taylor approximations of the Lagrangian which are formulated as Quadratic Programs (QP). Line-search or other methods are used to determine an appropriate step size. \blue{In theory, this solution method is well suited for being warm-started \cite{gill2012sequential}.} We will use KNITRO~\cite{byrd2006k} as the reference SQP solver.

\subsubsection{Interior-Point Methods (IPM)}
To deal with the constraint inequalities in the optimization problem, a logarithmic barrier term is added to the objective function with a multiplicative factor. This factor is decreased as the interior-point method converges, and the resulting barrier term resembles the indicator function. \blue{Interior point methods are challenging to efficiently warm start since the logarithmic barrier term initially keeps the solution away from inequality constraints that are binding at optimality \cite{forsgren2005warm}.} \bl{While theoretical performance guarantees on warm-starting interior-point methods for linear programs exist (see e.g. \cite{yildirim2002warm,benson2007exact}), common practice in the literature for non-convex programs (see, e.g., \cite{gould2002numerical}) is to follow an empirical approach based on extensive computational experiments.} We will use both KNITRO~\cite{byrd2006k} and IPOPT~\cite{wachter2006implementation} as reference IPM solvers \rev{as they are among the most robust and scalable solvers for AC-OPF \cite{kardos2018complete}}. 
\section{Simulations \& Results} 
\label{sec:results}
First, we specify the simulation setup. For the PGLib OPF test cases, we then evaluate the distance to AC feasibility and local optimality for the DC-OPF solution and solutions to the QC and SDP relaxations. We also study how these distances are correlated with the optimality gap. We use the DC-OPF as a computationally inexpensive benchmark. We next investigate the robustness and potential sub-optimality of penalization methods. Finally, we focus on warm-starting of non-convex solvers with solutions of inexact convex relaxations. \vspace{-0.3cm}
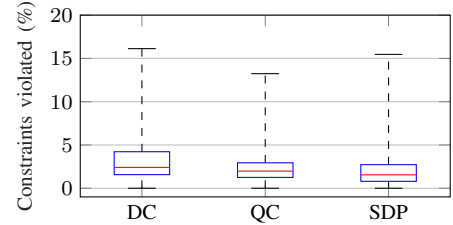
\begin{figure}
    \centering
    \begin{footnotesize}
%
%
\begin{tikzpicture}

\begin{axis}[%
width=6.5cm,
height=4cm,
xmin=0.5,
xmax=3.5,
xtick={1,2,3},
xticklabels={{DC},{QC},{SDP}},
ymin=-1,
ymax=20,
ylabel style={font=\color{white!15!black}},
ylabel={Constraints violated (\%)},
axis background/.style={fill=white},
ymajorgrids
]
\addplot [color=black, dashed, forget plot]
  table[row sep=crcr]{%
1	4.21979437844146\\
1	16.1290322580645\\
};
\addplot [color=black, dashed, forget plot]
  table[row sep=crcr]{%
2	2.93989251991598\\
2	13.2372214941022\\
};
\addplot [color=black, dashed, forget plot]
  table[row sep=crcr]{%
3	2.72375544552876\\
3	15.4652686762779\\
};
\addplot [color=black, dashed, forget plot]
  table[row sep=crcr]{%
1	0\\
1	1.57296625329735\\
};
\addplot [color=black, dashed, forget plot]
  table[row sep=crcr]{%
2	0\\
2	1.25422668829035\\
};
\addplot [color=black, dashed, forget plot]
  table[row sep=crcr]{%
3	0\\
3	0.802847414976811\\
};
\addplot [color=black, forget plot]
  table[row sep=crcr]{%
0.8875	16.1290322580645\\
1.1125	16.1290322580645\\
};
\addplot [color=black, forget plot]
  table[row sep=crcr]{%
1.8875	13.2372214941022\\
2.1125	13.2372214941022\\
};
\addplot [color=black, forget plot]
  table[row sep=crcr]{%
2.8875	15.4652686762779\\
3.1125	15.4652686762779\\
};
\addplot [color=black, forget plot]
  table[row sep=crcr]{%
0.8875	0\\
1.1125	0\\
};
\addplot [color=black, forget plot]
  table[row sep=crcr]{%
1.8875	0\\
2.1125	0\\
};
\addplot [color=black, forget plot]
  table[row sep=crcr]{%
2.8875	0\\
3.1125	0\\
};
\addplot [color=blue, forget plot]
  table[row sep=crcr]{%
0.775	1.57296625329735\\
0.775	4.21979437844146\\
1.225	4.21979437844146\\
1.225	1.57296625329735\\
0.775	1.57296625329735\\
};
\addplot [color=blue, forget plot]
  table[row sep=crcr]{%
1.775	1.25422668829035\\
1.775	2.93989251991598\\
2.225	2.93989251991598\\
2.225	1.25422668829035\\
1.775	1.25422668829035\\
};
\addplot [color=blue, forget plot]
  table[row sep=crcr]{%
2.775	0.802847414976811\\
2.775	2.72375544552876\\
3.225	2.72375544552876\\
3.225	0.802847414976811\\
2.775	0.802847414976811\\
};
\addplot [color=red, forget plot]
  table[row sep=crcr]{%
0.775	2.4076921278163\\
1.225	2.4076921278163\\
};
\addplot [color=red, forget plot]
  table[row sep=crcr]{%
1.775	1.97955259865108\\
2.225	1.97955259865108\\
};
\addplot [color=red, forget plot]
  table[row sep=crcr]{%
2.775	1.54917110917608\\
3.225	1.54917110917608\\
};
\end{axis}
\end{tikzpicture}%
    \end{footnotesize}
    \vspace{-0.3cm}
    \caption{The overall share of constraints violated in the AC-PF with the generators' setpoints for active power outputs and voltage magnitudes fixed to the values obtained from the DC-OPF, the QC relaxation, and the SDP relaxation \bl{for the successfully convergent 84 instances.}}
        \vspace{-0.4cm}
    \label{Mean_Max_Viol}
\end{figure}

\subsection{Simulation Setup}
\begin{figure*}[!th]
    \begin{footnotesize}
   \input{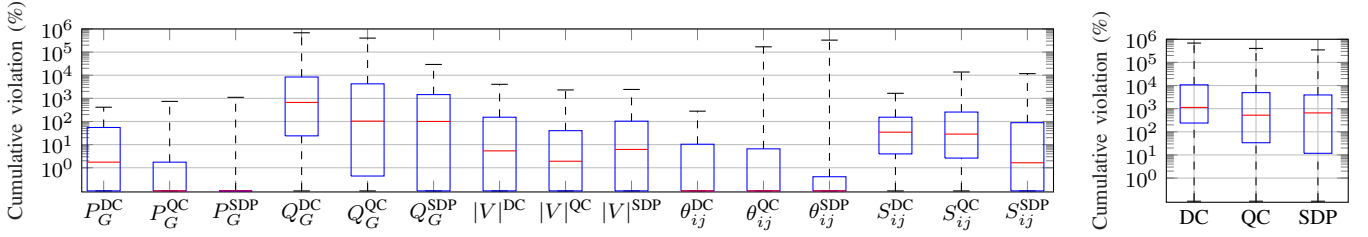} \quad  
%
%
\begin{tikzpicture}

\begin{axis}[%
width=4cm,
height=3.75cm,
xmin=0.55184331797235,
xmax=3.48963133640553,
xtick={1,2,3},
xticklabels={{DC},{QC},{SDP}},
ymin=0.09,
ymax=1000000,
ymode=log,
ytick = {1,10,100,1000,10000,100000,1000000},
ylabel style={font=\color{white!15!black}},
ylabel={Cumulative violation (\%)},
axis background/.style={fill=white},
xmajorgrids,
ymajorgrids
]
\addplot [color=black, dashed, forget plot]
  table[row sep=crcr]{%
1	10742.4231578479\\
1	685712.305897626\\
};
\addplot [color=black, dashed, forget plot]
  table[row sep=crcr]{%
2	4931.51178684263\\
2	403290.011869026\\
};
\addplot [color=black, dashed, forget plot]
  table[row sep=crcr]{%
3	3918.66761420272\\
3	347868.45638444\\
};
\addplot [color=black, dashed, forget plot]
  table[row sep=crcr]{%
1	0.1\\
1	239.459847974271\\
};
\addplot [color=black, dashed, forget plot]
  table[row sep=crcr]{%
2	0.1\\
2	33.6720109975326\\
};
\addplot [color=black, dashed, forget plot]
  table[row sep=crcr]{%
3	0.1\\
3	11.772593302893\\
};
\addplot [color=black, forget plot]
  table[row sep=crcr]{%
0.8875	685712.305897626\\
1.1125	685712.305897626\\
};
\addplot [color=black, forget plot]
  table[row sep=crcr]{%
1.8875	403290.011869026\\
2.1125	403290.011869026\\
};
\addplot [color=black, forget plot]
  table[row sep=crcr]{%
2.8875	347868.45638444\\
3.1125	347868.45638444\\
};
\addplot [color=black, forget plot]
  table[row sep=crcr]{%
0.8875	0.1\\
1.1125	0.1\\
};
\addplot [color=black, forget plot]
  table[row sep=crcr]{%
1.8875	0.1\\
2.1125	0.1\\
};
\addplot [color=black, forget plot]
  table[row sep=crcr]{%
2.8875	0.1\\
3.1125	0.1\\
};
\addplot [color=blue, forget plot]
  table[row sep=crcr]{%
0.775	239.459847974271\\
0.775	10742.4231578479\\
1.225	10742.4231578479\\
1.225	239.459847974271\\
0.775	239.459847974271\\
};
\addplot [color=blue, forget plot]
  table[row sep=crcr]{%
1.775	33.6720109975326\\
1.775	4931.51178684263\\
2.225	4931.51178684263\\
2.225	33.6720109975326\\
1.775	33.6720109975326\\
};
\addplot [color=blue, forget plot]
  table[row sep=crcr]{%
2.775	11.772593302893\\
2.775	3918.66761420272\\
3.225	3918.66761420272\\
3.225	11.772593302893\\
2.775	11.772593302893\\
};
\addplot [color=red, forget plot]
  table[row sep=crcr]{%
0.775	1132.75991912228\\
1.225	1132.75991912228\\
};
\addplot [color=red, forget plot]
  table[row sep=crcr]{%
1.775	524.039117462449\\
2.225	524.039117462449\\
};
\addplot [color=red, forget plot]
  table[row sep=crcr]{%
2.775	655.744971290321\\
3.225	655.744971290321\\
};
\end{axis}
\end{tikzpicture}%
    \end{footnotesize}
       \vspace{-0.3cm}
    \caption{\bl{Cumulative normalized constraint violation for the AC-PF solution with the generators' setpoints for active power outputs and voltage magnitudes fixed to the values obtained from the DC-OPF, the QC relaxation, and the SDP relaxation, considering constraints corresponding to each state variable  $P_G$, $Q_G$, $|V|$, $\Theta_{ij}$, $S_{ij}$ (left figure), and accumulated for all state variables (right figure) for the successfully convergent 84 instances.}}
    \label{Mean_Viol_all}
    \vspace{-0.2cm}
\end{figure*}
\begin{figure*}[!th]
    \centering
    \begin{footnotesize}
   \input{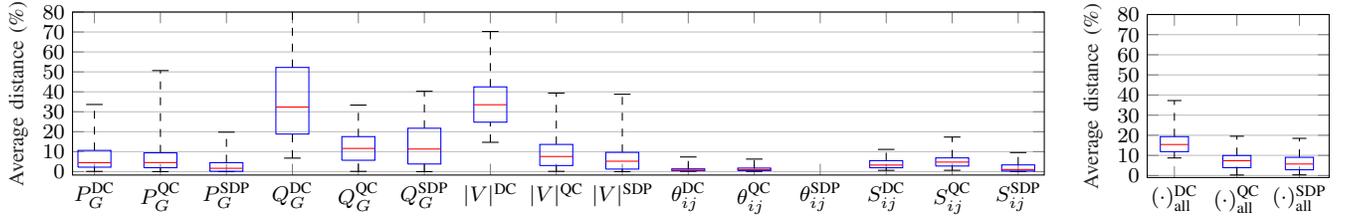}  \quad  
%
%
\begin{tikzpicture}

\begin{axis}[%
width=4cm,
height=3.75cm,
xmin=0.565668202764977,
xmax=3.4758064516129,
xtick={1,2,3},
xticklabels={$(\cdot)^{\text{DC}}_{\text{all}}$,$(\cdot)^{\text{QC}}_{\text{all}}$,$(\cdot)^{\text{SDP}}_{\text{all}}$},
ymin=-1,
ymax=80,
ymajorgrids,
ytick={0, 10, 20, 30, 40, 50, 60, 70, 80},
ylabel style={font=\color{white!15!black}},
ylabel={Average distance (\%)},
axis background/.style={fill=white}
]
\addplot [color=black, dashed, forget plot]
  table[row sep=crcr]{%
1	19.3310784963725\\
1	37.2471077026277\\
};
\addplot [color=black, dashed, forget plot]
  table[row sep=crcr]{%
2	9.94679543053748\\
2	19.5244638367425\\
};
\addplot [color=black, dashed, forget plot]
  table[row sep=crcr]{%
3	9.04364550384405\\
3	18.4756239598314\\
};
\addplot [color=black, dashed, forget plot]
  table[row sep=crcr]{%
1	8.73678028549369\\
1	11.859206821313\\
};
\addplot [color=black, dashed, forget plot]
  table[row sep=crcr]{%
2	0.258115811190433\\
2	3.8927902657898\\
};
\addplot [color=black, dashed, forget plot]
  table[row sep=crcr]{%
3	0.326291061459929\\
3	2.85549046953624\\
};
\addplot [color=black, forget plot]
  table[row sep=crcr]{%
0.8875	37.2471077026277\\
1.1125	37.2471077026277\\
};
\addplot [color=black, forget plot]
  table[row sep=crcr]{%
1.8875	19.5244638367425\\
2.1125	19.5244638367425\\
};
\addplot [color=black, forget plot]
  table[row sep=crcr]{%
2.8875	18.4756239598314\\
3.1125	18.4756239598314\\
};
\addplot [color=black, forget plot]
  table[row sep=crcr]{%
0.8875	8.73678028549369\\
1.1125	8.73678028549369\\
};
\addplot [color=black, forget plot]
  table[row sep=crcr]{%
1.8875	0.258115811190433\\
2.1125	0.258115811190433\\
};
\addplot [color=black, forget plot]
  table[row sep=crcr]{%
2.8875	0.326291061459929\\
3.1125	0.326291061459929\\
};
\addplot [color=blue, forget plot]
  table[row sep=crcr]{%
0.775	11.859206821313\\
0.775	19.3310784963725\\
1.225	19.3310784963725\\
1.225	11.859206821313\\
0.775	11.859206821313\\
};
\addplot [color=blue, forget plot]
  table[row sep=crcr]{%
1.775	3.8927902657898\\
1.775	9.94679543053748\\
2.225	9.94679543053748\\
2.225	3.8927902657898\\
1.775	3.8927902657898\\
};
\addplot [color=blue, forget plot]
  table[row sep=crcr]{%
2.775	2.85549046953624\\
2.775	9.04364550384405\\
3.225	9.04364550384405\\
3.225	2.85549046953624\\
2.775	2.85549046953624\\
};
\addplot [color=red, forget plot]
  table[row sep=crcr]{%
0.775	15.4040216479841\\
1.225	15.4040216479841\\
};
\addplot [color=red, forget plot]
  table[row sep=crcr]{%
1.775	7.33636307303388\\
2.225	7.33636307303388\\
};
\addplot [color=red, forget plot]
  table[row sep=crcr]{%
2.775	5.72134371610222\\
3.225	5.72134371610222\\
};
\end{axis}
\end{tikzpicture}%
    \end{footnotesize}
       \vspace{-0.3cm}
    \caption{Averaged normalized distance of the solution variables of the DC-OPF, QC and SDP relaxations to the locally optimal solution reported by IPOPT considering each state variable  $P_G$, $Q_G$, $|V|$, $\Theta_{ij}$, $S_{ij}$ (left figure), and averaged for all state variables (right figure) \bl{for all considered 96 test cases}. \blueEPSR{Note that the voltage angles are not available from \mbox{PowerModels.jl} if the SDP relaxation is not exact. The underlying reason is that the voltage angles cannot be directly extracted if the matrix variable $W$ is not rank-1.}}
    \label{Mean_Dist}
    \vspace{-0.55cm}
\end{figure*}

We use the implementations of the AC-OPF, DC-OPF, QC relaxation, and SDP relaxation provided in \mbox{PowerModels.jl}~\cite{coffrin2018powermodels}, a computationally efficient open-source implementation in Julia. In \mbox{PowerModels.jl}, we use KNITRO and IPOPT to solve the non-convex AC-OPF, MOSEK to solve the DC-OPF and the SDP relaxation, and IPOPT to solve the QC relaxation. The analysis in this work uses the PGLib OPF Benchmarks v18.08~\cite{PGLIB}, in particular, the test cases ranging from 14 to 3120 buses under typical, congested, and small angle difference conditions. We exclude the test case \textit{case2000\_tamu} since the SDP relaxation fails for this test case. For the remaining 96 test cases, the QC and SDP relaxations return the optimal solution, although it should be noted that MOSEK reports ``stall'' in some of the test cases due to numerical issues. For the small angle difference conditions, the DC-OPF is infeasible for several instances. For these, we iteratively relax the angle difference constraints in the DC-OPF problem in 10\% steps until we obtain a feasible solution. All simulations are carried out on a laptop.
\subsection{Distances to AC Feasibility and Local Optimality}

\subsubsection{AC Feasibility}

\label{AC-Feasbility} 
\blue{In the following, we evaluate the constraint violation resulting from AC power flow solutions in M{\sc atpower}~\cite{zimmerman2011matpower} obtained using the generators' setpoints for active power outputs and voltage magnitudes from the DC-OPF, QC relaxation, and SDP relaxation. \blueEPSR{Note that the generator voltage set-points are fixed to 1.0~p.u. for the DC-OPF.} For this purpose, we make the following assumptions: The largest generator is selected as slack bus. A numerical tolerance of 0.1\% is considered as minimum constraint violation limit. \blueEPSR{Note that out of 96 test cases considered, the AC power flow does not converge for four test cases under any of the three loading configurations due to numerical ill-conditioning, i.e., in total 12 AC power flows do not converge.} The characteristics of these four test cases, which model parts of the French transmission network, are detailed in~\cite{josz2016ac}. These AC power flows also do not converge using PowerModels.jl with IPOPT. The underlying reasons for power flow non-convergence are outside the scope of this work. The following AC feasibility analysis uses the results of the M{\sc atpower} AC-PF and focuses on the remaining 84 convergent test cases.}

\blue{In Fig.~\ref{Mean_Max_Viol}, we investigate the fraction of constraints violated and show that, in terms of fraction of constraints violated, the lower 75th percentiles for the QC and SDP relaxations (2.9\% and 2.7\%, respectively) are lower than for DC-OPF (4.2\%). \bl{In Fig.~\ref{Mean_Viol_all}, we show the cumulative normalized constraint violation for each constraint type and for the entire AC-OPF problem. We can make the following observations: First, for active power constraints, the 75th percentile for the SDP relaxation do not exhibit constraint violation, whereas for the DC-OPF this holds true only for the 25th percentile. The cumulative violation of generator reactive power constraints represents a significant share of the overall cumulative violation, with the SDP and QC relaxations exhibiting lower constraint violations than DC-OPF. For the voltage magnitude, all three initializations of the AC-PF lead to similar cumulative violations. For the apparent branch flow constraints, it can be observed that AC power flows initialized with the solution to the SDP relaxation lead to significantly smaller constraint violations compared to DC-OPF and QC relaxation.} Regarding the overall cumulative constraint violation, the lower 25th percentile for the SDP relaxation (11.8\%) is lower than the QC relaxation (33.7\%) and the DC-OPF (239.5\%), as these test cases are very close to a rank-1 solution. Furthermore, the lower 75th percentile of the cumulative constraint violation for the QC and SDP relaxations are 54.1\% and 63.5\% lower than for the DC-OPF, respectively, highlighting that the obtained solutions are closer to AC feasibility. We also evaluate the type of constraint leading to the maximum constraint violation and we find that this is the generator reactive power in 71.1\%, 50.6\% and 49.4\% of test cases for the DC-OPF, QC relaxation and SDP relaxation, respectively. If we enforce the reactive power limits in the AC-PF, whenever possible, active generator power and voltage limit violations occur more often instead. \bl{Developing systematic procedures to convert $PV$ to $PQ$ buses when generator reactive power limits are violated in the AC power flow is an interesting direction for future work.} The cumulative constraint violation relates to the distance to an AC-feasible solution, which is not necessarily the same as the distance to a local optimum. We analyze the latter next.}
\subsubsection{Local Optimality}
For the 96 considered test cases, we compute the average normalized distance between the locally optimal solution found by the non-convex solver IPOPT and the solutions to the DC-OPF, the QC relaxation, and the SDP relaxation. Fig.~\ref{Mean_Dist} shows the average distance for each state variable and for the average of all state variables. Considering all state variables, the lower 75th percentiles for the SDP and QC relaxations are less than 10\% (9.0\% and 9.9\%), and significantly smaller than the DC-OPF solution (19.3\%). Looking at the individual state variables, the SDP and QC relaxations' solutions are significantly closer to the local solution than the DC-OPF, particularly for the reactive generator power $Q_G$ and voltage magnitude $|V|$. For generator active power and the apparent branch flow, the SDP relaxation is the closest (less than 5\% for the lower 75th percentile). 
Since obtaining the voltage angles from the SDP relaxation's solutions is not straightforward, the angles are set to zero, and we do not report the distance in Fig.~\ref{Mean_Dist}. The closeness of solutions from the QC and SDP relaxations relative to the local solution motivates our investigation of warm-starting techniques in Section~\ref{Warm_Start}.
\subsubsection{Correlation with Optimality Gap}
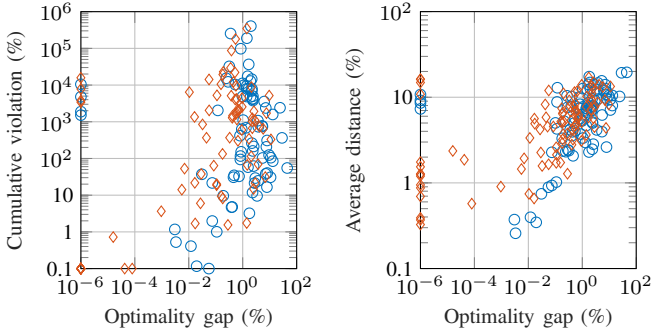
\begin{figure}
    \centering
    \begin{footnotesize}
%
%
\definecolor{mycolor1}{rgb}{0.00000,0.44700,0.74100}%
\definecolor{mycolor2}{rgb}{0.85000,0.32500,0.09800}%
\begin{tikzpicture}

\begin{axis}[%
width=4.45cm,
height=5cm,
xmode=log,
xmin=1e-06,
xmax=100,
xtick = {1e-06, 1e-04, 1e-02, 1e00, 1e02},
xlabel style={font=\color{white!15!black}},
xlabel={Optimality gap (\%)},
ymin=0.1,
ymax=1000000,
ytick = {0.1, 1, 10, 100,1000,10000,100000,1000000},
yticklabels={0.1, 1, 10, $10^2$,$10^3$,$10^4$,$10^5$,$10^6$},
ymode = log,
ylabel style={font=\color{white!15!black}},
ylabel={Cumulative violation (\%)},
axis background/.style={fill=white},
xminorticks=true,
xminorgrids = true,
xmajorgrids = true,
xminorgrids = true,
ymajorgrids = true,
legend style={legend pos = north west, legend cell align=left, align=left, draw=white!15!black}
]
\addplot [color=mycolor1, draw=none, mark=o, mark options={solid, mycolor1}]
  table[row sep=crcr]{%
0.109078537812735	1.00565453908936\\
0.0120340035907329	0.40854666169136\\
0.055670455720358	0.1\\
0.386009961925471	4.70366528497972\\
10.7785646761365	30.7585714344847\\
0.480564624406921	34.7219578333804\\
0.455152464005582	30.8090124080564\\
0.0305131507658318	37.1109124032255\\
0.730476600635621	32.6220641616847\\
2.19009979552861	261.681014989186\\
7.53983139358844	207.39320511549\\
0.119502671685545	2041.41835355696\\
0.00335014157012425	0.526648910107442\\
3.80369351791656	730.416600430354\\
2.55386483398785	212.712335816195\\
5.38048711888349	742.752078450615\\
1.67395000554676	8195.22403629162\\
2.39149599135425	3791.57866408687\\
2.21135372336195	7792.14770912126\\
0.989623402296724	11141.1406642005\\
0.287933842860932	12239.8075961217\\
0.245048381476742	15284.1246738168\\
0.355387078381275	253525.34661659\\
0.317449287658245	32335.8724548617\\
1.42742298810131	189193.553165593\\
1.06506792875267	9797.52381412574\\
0.978021487393987	4890.61103014687\\
0.532909339883203	10314.8391283793\\
1.76801515406821	3.21194120583417\\
13.0052838978833	84.3338135989186\\
44.6009306644326	55.0404063360838\\
2.75440256655376	21.4212371262191\\
3.72403762600474	24.5632634065669\\
1.56919250196695	58.0134387896926\\
0.0743522084170167	21.5467709370796\\
11.0615526210842	105.077427016687\\
8.12653333452742	40.1643289056482\\
28.6267563697924	364.037041240652\\
5.43793440214143	103.247369070247\\
7.17452375762347	2378.73159064921\\
0.0195941965142055	0.116294092058302\\
0.792351185219442	418.69719550192\\
0.875335636921482	152.170106912248\\
0.0742313509582648	1.97835438388157\\
0.9246761922034	629.381039422978\\
0.846428135423893	97722.4752238486\\
1.9520263939936	4760.26230940751\\
1e-06	1515.58602343129\\
12.9667775608613	1581.4713169531\\
3.66244224368891	1073.54167693612\\
1e-06	10505.6873289734\\
1e-06	4788.3143778492\\
2.30352699888131	4258.54858332342\\
1.31131702808683	683.257358204753\\
1e-06	1882.56853188698\\
24.1419865650632	2408.27883877379\\
7.15710074326212	38.1687139554506\\
2.92704451303031	29.251037026267\\
2.31339577853624	11.8261856704229\\
0.407559419575632	4.78493234607197\\
3.41015716487523	9.56971084603154\\
0.19001070276109	9.62687416798105\\
0.828492189730579	116.560158411759\\
2.53419739648274	61.9570328862395\\
0.817550451278348	94.0108744148792\\
9.47915864084228	124.117605776837\\
8.01663314506471	309.693490867971\\
1.04110134555659	1110.02348506496\\
0.00304074724752912	1.16804354558015\\
5.23507062302724	772.849391566751\\
2.35357538786833	160.519553591981\\
7.89467546360989	16.5983590621098\\
6.25526319718693	3554.95767588153\\
2.36214024873803	4094.63774828498\\
2.17243615326042	6916.56649956069\\
2.15117946820135	14527.1965353489\\
1.52704508302259	17390.4249656122\\
1.92027508946965	12897.2226629952\\
1.99971451990063	403290.011869026\\
1.67190764996117	36600.4922684616\\
2.38916972626324	4972.41254353838\\
1.40515522877547	9490.73243209654\\
1.40170621033998	4867.99441362025\\
1.41489940787125	7921.15084993636\\
};

\addplot [color=mycolor2, draw=none, mark=diamond, mark options={solid, mycolor2}]
  table[row sep=crcr]{%
1e-06	0.1\\
1e-06	0.1\\
1e-06	0.1\\
1e-06	0.1\\
7.98612957431111e-05	0.1\\
0.00553498283579756	14.1509949246043\\
1e-06	0.1\\
1e-06	0.1\\
0.00696962431802062	52.8953951079177\\
0.178577981217332	133.166580664372\\
2.25999479541197	1017.76425266157\\
0.057015482446876	2044.59595275669\\
1e-06	0.1\\
2.27045926778194	668.59296755911\\
0.108501722856125	19.9871483333545\\
2.10770360765428	8810.37916051482\\
0.665547451438131	303.37651362664\\
0.555758648701354	1616.50036860033\\
0.728151900711571	2102.67884631593\\
0.375435318529405	5119.34467057797\\
0.016349615118505	1355.095073401\\
0.0324451218437383	848.07694031674\\
0.0564005944391877	14216.5443454105\\
0.0105413570734014	6414.68147645966\\
0.542573800116042	179675.62656675\\
0.410494727881217	1283.16073260017\\
0.156596191308722	10550.1427821953\\
0.111203762330869	5250.48461694086\\
1e-06	0.1\\
2.06463816119666	114.019012508119\\
1.40918988222286	1.74161684455743\\
0.283201444022962	1.54933271519599\\
1.58145686213196e-05	0.721737868478874\\
0.181800782705765	8.17775792743023\\
0.000956214783520704	3.64939090033758\\
2.91209720937635	21.4061270040456\\
6.88182626672409	31.1035858308369\\
11.1384901802208	330.291999345366\\
1.70970347045727	343.138600825278\\
0.641437870862749	2122.28623357271\\
1e-06	0.1\\
0.328546799393026	368.461480630602\\
0.0290109876383116	37.8522704779548\\
1e-06	0.1\\
0.528149522039212	4051.14402551121\\
0.381979846458091	87148.9511214807\\
0.608148536918596	1155.63383751442\\
1e-06	3425.24676325557\\
2.60207258358697	1925.80190593475\\
2.83291751004656	1100.73202762937\\
1e-06	15801.3529505114\\
1e-06	7815.90003575023\\
0.945271363453082	4408.7882600166\\
0.921794453483915	597.792231293806\\
1e-06	3786.19120289423\\
9.58528111077214	5164.92422100358\\
0.0272015814298143	6.0006530599757\\
2.51999368497783	14.7914987438944\\
0.160938599809279	9.39419168118169\\
0.0172971777243625	1.69642924753037\\
4.21067457123669e-05	0.1\\
0.0167492325199636	21.7120258376649\\
0.0845997050993641	63.8414368134924\\
1.47198820097605	66.8989648517405\\
0.0453826688804226	367.057926549237\\
3.7009009103815	69.9924034546691\\
2.27612660065885	1246.3218150897\\
0.951217556385009	2558.3476275313\\
1e-06	0.1\\
4.16790753337407	732.585063094585\\
0.108008574094431	18.7249735036661\\
7.58713025886821	1718.95523117156\\
5.75127490575136	642.896975021531\\
0.589158768103659	1758.45237443449\\
0.729396096163626	2103.09832476031\\
0.560155291747166	9594.01506525969\\
0.183999563559534	19389.4934546275\\
0.514394561910325	22193.1711133662\\
0.356225808264921	29368.2624232061\\
0.192102663613924	22986.3861326947\\
1.41806412876284	347868.45638444\\
0.421313343633067	1458.88233362371\\
0.450569923461352	12133.1007988854\\
0.603306393870529	4202.03798916257\\
};

\end{axis}
\end{tikzpicture}%
 \, 
%
%
\definecolor{mycolor1}{rgb}{0.00000,0.44700,0.74100}%
\definecolor{mycolor2}{rgb}{0.85000,0.32500,0.09800}%
\begin{tikzpicture}

\begin{axis}[%
width=4.45cm,
height=5cm,
xmode=log,
xmin=1e-06,
xmax=100,
xtick = {1e-06, 1e-04, 1e-02, 1e00, 1e02},
xlabel style={font=\color{white!15!black}},
xlabel={Optimality gap (\%)},
ymin=0.1,
ymax=100,
ytick = {0.1, 1, 10, 100},
yticklabels={0.1, 1, 10, $10^2$},
ymode=log,
ylabel style={font=\color{white!15!black}},
ylabel={Average distance (\%)},
axis background/.style={fill=white},
xminorticks=true,
xminorgrids = true,
xmajorgrids = true,
xminorgrids = true,
ymajorgrids = true,
scaled y ticks=false,
legend style={legend pos = north west, legend cell align=left, align=left, draw=white!15!black}
]
\addplot [color=mycolor1, draw=none, mark=o, mark options={solid, mycolor1}]
  table[row sep=crcr]{%
0.109078537812735	1.02174423863648\\
0.0120340035907329	0.396787452564712\\
0.055670455720358	0.890377334723706\\
0.386009961925471	7.82189752457966\\
10.7785646761365	9.39113590031756\\
0.480564624406921	4.96129797299161\\
0.455152464005582	2.07131243302707\\
0.0305131507658318	0.745198650817912\\
0.730476600635621	7.45523838534879\\
2.19009979552861	7.02993622560367\\
7.53983139358844	10.8193841679925\\
0.119502671685545	12.9111177601128\\
0.00335014157012425	0.258115811190433\\
3.80369351791656	14.6715247381257\\
2.55386483398785	2.83118054753895\\
5.38048711888349	3.9655702160603\\
1.67395000554676	8.61740080877595\\
2.39149599135425	4.99851721002041\\
1.81118251811911	17.4100132508903\\
0.114026069449003	4.27456235042897\\
2.21135372336195	7.04559605329256\\
0.989623402296724	6.24469543243165\\
0.287933842860932	2.43392553713752\\
0.245048381476742	2.31377170954731\\
0.355387078381275	4.18466224007214\\
0.317449287658245	2.37932775709959\\
0.115652809681488	8.81069659919232\\
1.42742298810131	7.3495417350073\\
0.101009724065504	5.01369647019325\\
1.06506792875267	3.33541703540283\\
0.978021487393987	5.28491306941936\\
0.532909339883203	4.03426882244279\\
1.76801515406821	2.29395135201166\\
13.0052838978833	10.484317191552\\
44.6009306644326	19.5244638367425\\
2.75440256655376	11.8303198066111\\
3.72403762600474	9.51642865919347\\
1.56919250196695	7.98737398934477\\
0.0743522084170167	2.27674322326817\\
11.0615526210842	7.09652342275984\\
8.12653333452742	11.569560748592\\
28.6267563697924	19.2988536159544\\
5.43793440214143	6.20196985162746\\
7.17452375762347	12.6264467401025\\
0.0195941965142055	0.346148531162441\\
0.792351185219442	10.0638192569549\\
0.875335636921482	3.71505437113071\\
0.0742313509582648	0.9334551198908\\
0.9246761922034	5.39792413748264\\
0.846428135423893	5.64385609345438\\
0.281104601098514	7.6048575540383\\
0.432087060638009	3.94913638871486\\
1.9520263939936	7.4980225055429\\
1e-06	8.56711255915934\\
12.9667775608613	8.00518786010888\\
3.66244224368891	12.1484376745263\\
1e-06	10.7386612507022\\
1e-06	9.18307170650733\\
0.210263397836719	5.25739097556587\\
2.30352699888131	7.88139524135897\\
0.179846512572612	3.31294407117804\\
1.31131702808683	2.96190376593138\\
1e-06	7.30729822126058\\
24.1419865650632	10.2179970779853\\
7.15710074326212	14.0183112895318\\
2.92704451303031	9.17853487429662\\
2.31339577853624	9.59520346108254\\
0.407559419575632	10.1944935259739\\
3.41015716487523	4.15667019619516\\
0.19001070276109	2.35848044623459\\
0.828492189730579	7.44232840171105\\
2.53419739648274	7.91627473315535\\
0.817550451278348	8.28237901716269\\
9.47915864084228	10.3868886184151\\
8.01663314506471	9.7084863723596\\
1.04110134555659	12.6164240332245\\
0.00304074724752912	0.372616352994496\\
5.23507062302724	15.5380741929272\\
2.35357538786833	2.47232324018117\\
7.89467546360989	2.61857821065836\\
6.25526319718693	12.0931654037931\\
2.36214024873803	5.13897163722136\\
2.72321155771984	14.1581643187247\\
0.405396988278184	3.55174668746292\\
2.17243615326042	6.99036720172414\\
2.15117946820135	7.56244856828531\\
1.52704508302259	12.246461992522\\
1.92027508946965	14.8770357735208\\
1.99971451990063	12.568745930428\\
1.67190764996117	14.1192173675371\\
0.270935816007167	10.3993465973762\\
2.38916972626324	7.87707449521237\\
0.538028835365212	5.65607789028256\\
1.40515522877547	4.66934861219516\\
1.40170621033998	6.68884378542692\\
1.41489940787125	7.46800757777778\\
};

\addplot [color=mycolor2, draw=none, mark=diamond, mark options={solid, mycolor2}]
  table[row sep=crcr]{%
1e-06	0.946780853502886\\
1e-06	1.2718131674634\\
1e-06	1.20068832538034\\
1e-06	0.694800755648144\\
7.98612956987022e-05	0.572012860286458\\
0.00553498283575315	1.4509383159868\\
1e-06	0.874210066382455\\
1e-06	1.26095676439673\\
0.00696962431797621	1.93314155290518\\
0.178577981217276	5.0685822420189\\
2.25999479541192	9.28231792253205\\
0.0570154824468205	11.995147805832\\
1e-06	0.326291061459929\\
2.2704592677819	13.3827174695233\\
0.108501722856091	1.7154449410799\\
2.10770360765424	7.19349708402461\\
0.665547451438087	5.61056690831209\\
0.555758648701299	4.20683134981049\\
1.74172717360157	17.0016801461005\\
0.0113610146591392	3.34788118137674\\
0.728151900711527	5.82596472122279\\
0.37543531852936	6.76351040325997\\
0.0163496151184606	0.657287314893083\\
0.0324451218436939	1.57715055818862\\
0.056400594439121	1.64845944513896\\
0.010541357073357	0.744233963890919\\
0.0422432762707259	8.24289346156424\\
0.542573800115997	6.72969247880355\\
0.204559968772688	5.18838572765722\\
0.410494727881161	2.78994094788606\\
0.156596191308667	1.72807331909101\\
0.111203762330825	2.12537952647914\\
1e-06	1.75149545686726\\
2.06463816119663	9.70726714028213\\
1.40918988222281	3.63228127469662\\
0.283201444022918	1.2796327344464\\
1.58145685769107e-05	2.36090086527464\\
0.181800782705721	2.92103999118641\\
0.000956214783476295	0.904607908534658\\
2.91209720937632	7.7750797860218\\
6.88182626672404	9.28137013466403\\
11.1384901802207	13.3712718983802\\
1.7097034704572	5.81380962922429\\
0.641437870862704	11.7307931074529\\
1e-06	0.369232129670339\\
0.32854679939297	9.1791812913436\\
0.0290109876382671	2.26217774669618\\
1e-06	0.39053992709611\\
0.528149522039167	4.62488342489496\\
0.381979846458036	5.35824109696118\\
0.171953758501864	7.44447910998476\\
0.193806377257555	3.71281639543296\\
0.608148536918562	4.78213729297642\\
1e-06	16.1653659351433\\
2.60207258358691	5.67620933640732\\
2.83291751004652	4.48135421304781\\
1e-06	16.3071546998937\\
1e-06	14.8555378272614\\
0.0515923572480648	4.53866766684991\\
0.945271363453026	5.11571108130637\\
0.185911489711932	3.2163678161957\\
0.921794453483871	3.08396952069801\\
1e-06	10.8224791407213\\
9.5852811107721	8.23700171500509\\
0.0272015814297588	4.21079748425145\\
2.51999368497779	12.8444952138971\\
0.160938599809235	5.18552048314124\\
0.0172971777243069	4.91863792268367\\
4.2106745667958e-05	1.86139434493895\\
0.0167492325199081	5.62887780298015\\
0.0845997050993086	4.5564000025306\\
1.471988200976	10.2262823227918\\
0.0453826688803782	5.83104340597764\\
3.70090091038145	8.70714007783121\\
2.2761266006588	11.1145576960742\\
0.951217556384953	13.9637077428659\\
1e-06	1.54034317797966\\
4.16790753337404	13.9757129677167\\
0.108008574094398	2.60322812850589\\
7.58713025886818	4.38420730772238\\
5.75127490575131	13.1805611707228\\
0.589158768103604	4.86942294714992\\
2.66696082546163	14.9623261882797\\
0.280076082147951	3.77245175153878\\
0.729396096163581	6.25333284761062\\
0.560155291747111	8.46084789727969\\
0.1839995635595	6.41518051098023\\
0.514394561910281	11.4926074327988\\
0.356225808264854	7.86856778874662\\
0.192102663613891	6.64635101700302\\
0.200701391420721	11.1085334819802\\
1.41806412876281	7.76207003057011\\
0.505573509745416	7.72369812643536\\
0.421313343633012	4.2405238106206\\
0.450569923461319	3.20720243573252\\
0.603306393870495	6.90447711558096\\
};
\end{axis}
\end{tikzpicture}%
    \end{footnotesize}
    \vspace{-0.75cm}
    \caption{Correlation of the optimality gap with the cumulative constraint violation and the averaged distance to local optimality for the QC and SDP relaxations. Note that the axes are logarithmic. The red diamonds correspond to the SDP relaxation, and the blue circles correspond to the QC relaxation. \bl{Note that for the first metric we consider the 84 convergent cases, and for the second metric we consider all 96 test instances.}}
    \label{Opt_Gap_Dist_Viol}
    \vspace{-0.55cm}
\end{figure}
We investigate the correlation between the optimality gaps and both the averaged normalized distance and cumulative constraint violation for the QC and SDP relaxations as shown in Fig.~\ref{Opt_Gap_Dist_Viol}. Note that both axes are on a logarithmic scale and the values are thresholded at $10^{-6}$. We use the locally optimal solution obtained from IPOPT as the best known feasible point to compute the optimality gap. Even for cases with optimality gaps that are less than 1\%, both the distances to local optimality and the cumulative constraint violations can still be substantial, suggesting that the optimality gap does not adequately capture the tightness of a relaxation in terms of the decision variable accuracy.

Furthermore, there is a group of test cases with non-negligible distances to local optimality and substantial constraint violations that nevertheless have optimality gaps which are almost zero. The outliers are the following four test cases: \textit{case2383wp\_k\_\_api}, \textit{pglib\_opf\_case2746wop\_k\_\_api}, \textit{case2746wp\_k\_\_api}, \textit{case3012wp\_k\_\_api}, which are some of the test cases representing the Polish grid under congested operation conditions (api). \blueEPSR{A possible explanation is that the objective functions for these test cases are very flat with respect to the change in the active generator dispatch, i.e., there are many generators with the same cost parameters for these four test cases. Note that all four test cases only have linear cost parameters (the quadratic and constant cost terms are zero). To test this hypothesis, we created variants of these four test cases where we assigned a random linear cost function to each generator. To this end, we took random samples between the minimum and maximum linear cost parameters of the original test case. When recomputing optimality gaps and distances to AC-feasiblity and local optimality for these test cases, we observed that the obtained optimality gaps are non-zero, thus supporting our hypothesis.}

For the correlation analysis, we use Spearman's rank and Pearson's correlation coefficient. The first coefficient relates to the (possibly non-linear) monotonicity of the optimality gap with the two metrics described in Section~\ref{Metrics}. For the Pearson correlation coefficient, we identify the strength of possible linear relationships between the base-10 logarithms of the optimality gaps and these metrics. \blue{For both coefficients, the obtained values can range between -1 and 1, where the minimum and maximum values correspond to perfect negative or positive correlation, and the value of 0 expresses that the quantities are uncorrelated in this statistical measure. With regard to both the cumulative constraint violation and distance to local optimality metrics, the resulting Spearman's rank and Pearson correlation coefficients are in an interval between 0.32 and 0.59, showing both metrics are not strongly correlated with the optimality gap.} \vspace{-0.3cm}
\subsection{Penalization Methods for SDP Relaxation}
\label{Pen_results}
The previous section shows that the SDP relaxation's solution is close to both AC feasibility and local optimality in various test cases. To drive the solution towards rank-1 and consequently AC feasibility, we have presented three penalty terms from the literature in Section~\ref{Theory_Pen}. In this section, we provide a detailed analysis of the robustness of these heuristic penalization methods for recovering an AC-feasible solution \bl{and show how the proposed metrics can be used to assess the quality of the solutions obtained from penalized SDP relaxations}. We also quantify the sub-optimality incurred by modifying the original objective function. First, we present a five-bus test case serving as an illustrative example where penalty terms can fail to recover both an AC-feasible solution and a near-globally or locally optimal solution. Then, we provide a detailed numerical analysis on the PGLib OPF test case database considering all test cases with up to 300 buses. To facilitate comparability, as the absolute objective function values of the different test cases vary significantly, we define the penalty weight in percent of the original objective function value $f_{\text{cost}}^0$ of the SDP relaxation with no penalty term included. As an illustrative example, a penalty weight of $\epsilon_{\text{pen}} = 1\%$ corresponds to $\epsilon_{\text{pen}} = 0.01 \times f_{\text{cost}}^0$. \bl{While previous works have used heuristic measures for the rank-1 property of the obtained $W$ matrix, in particular requiring the ratio of the first and second eigenvalue of the obtained matrix $W$ to be at least between $10^4$ and $10^6$, we directly compute the cumulative normalized constraint violation defined in \eqref{eq:constrviol} to assess whether the obtained decision variables are AC-feasible or not. This is based on our observations that despite having an eigenvalue ratio of larger than $10^4$, several of the test cases exhibit non-zero distances to AC-feasibility, i.e., constraint violations occur. Note that, due to numerical inaccuracies, if the cumulative normalized constraint violation is below $0.1$\%, we assume that no constraint violations occur, i.e., the solution is AC-feasible.}

\subsubsection{Five-Bus Test Case}
\begin{figure}
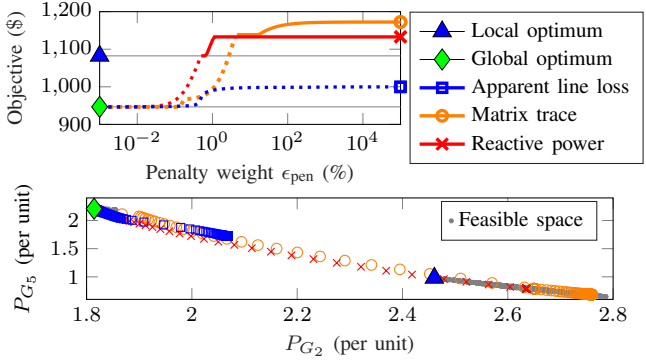

\begin{center}
		\begin{footnotesize}
			\input{./Figures/Pen_Comp_5bus.tex} \\
			\vspace{-0.2cm}
			\input{./Figures/5bus_FS.tex}  
		\end{footnotesize}
		\end{center}
		\vspace{-16pt}
	\caption{The upper plot shows a comparison of penalization methods applied to the SDP relaxation for the five-bus test case from~\cite{Bukhsh2013}. The solid lines represent the regions within which a rank-1 solution matrix is obtained and hence the SDP solution is AC-feasible. Conversely, the dotted lines correspond to higher-rank solution matrix for which the SDP solution is not AC-feasible. The two horizontal lines indicate the objective values of the local and global optima. The bottom plot shows the results from the penalization methods projected on the disconnected feasible space with respect to the active power generation $P_{G_2}$ and $P_{G_5}$. The feasible space is reproduced from~\cite{Molzahn2017_FS}. Please note that the trajectory associated with the matrix trace penalty term does not pass through the spurious locally optimal solution but leads to a different part of the feasible space.} 
	\label{5bus_Pen}
	\vspace{-0.4cm}
\end{figure}
We investigate the five-bus test case from~\cite{Bukhsh2013}. The feasible space of this system is visualized in~\cite{Molzahn2017_FS} and shown to be disconnected with one local solution in addition to the global optimum. The upper plot in Fig.~\ref{5bus_Pen} shows the objective value versus the penalty weight for the three different penalty functions. We use a \rev{fine} penalty step size $\Delta \epsilon_{\text{pen}}$ of $10^{0.05}$\% ranging from $10^{-3}$\% to $10^5$\%. The solid line sections represent the region within which each penalty term yields a rank-1 solution matrix $W$. The lower plot in Fig.~\ref{5bus_Pen} shows the feasible space and the corresponding results of the penalization methods projected onto the disconnected feasible space with respect to the active power generation $P_{G_1}$ and $P_{G_2}$. We make the following observations: The SDP relaxation without a penalty term included is inexact, but very close to the global optimum. Excluding the voltage angles (\blueEPSR{which are not available from the PowerModels.jl implementation of the SDP relaxation}), the normalized distances to the local and global optima are 4.2\% and 0.1\%, respectively. The M{\sc atpower} AC-PF, however, does not converge if initialized with the solution of the inexact SDP relaxation. The performance of the three penalization terms is compared: By including a reactive power penalty, the solution is moving towards the local optimum and away from the global optimum. If the reactive power penalty weight is chosen high enough, we obtain a rank-1 matrix, represented by the solid part of the line, and we obtain the locally optimal solution for a small interval of penalty weights ($\epsilon_{\text{pen}} \approx 0.5\%$ to $0.7\%$). As we increase the penalty weight beyond $0.7\%$, the solution is driven towards a different rank-1 (and hence AC-feasible) point which incurs sub-optimalities of 19.7\% and 4.7\% for the globally and locally optimal solutions, respectively. For the matrix trace, a sufficiently high penalty term of $4.5 \%$ yields a rank-1 solution. In this case, we cannot obtain the locally optimal solution, but incur at least a sub-optimality of $20.3$\% and $5.2$\% with respect to the global and local optima, respectively. As seen in the plot of the feasible space, increasingly penalizing the matrix trace results in movement toward the same portion of the feasible space as increasing the reactive power penalty, but does not result in passing through the locally optimal solution. The apparent branch loss penalty term fails to recover a rank-1 solution.

\subsubsection{PGLib OPF Test Cases up to 300 Buses}
\begin{table}[]
\caption{\bl{Performance of penalization methods on PGLib OPF test cases}}
    \centering
    \setlength\tabcolsep{3pt} 
\begin{tabular}{l c c c c c}
\toprule 

         Penalty   & Add. rank-1  & Range of & min. opt. & Range of & max. opt.\\
    term & cases (\%) &  $\epsilon_{\text{pen}}^{\text{min}}$ (\%)& gap (\%)&  $\epsilon_{\text{pen}}^{\text{max}}$ (\%)&  gap (\%)\\
     \midrule 
     $Q_G$& 42.2 & $10^{-5}$--$10^{9}$  & 0.0--24.0 &  $10^{8}$--$10^{10}$ &  0.2--39.2 \\
     $\text{Tr}\{W\}$ & 17.8 & $10^{-5}$--$10^{7}$   & 0.0-24.0 & $10^{5}$--$10^{10}$& 0.0-31.0\\
     $| S_{ij} - S_{ji}|$ & 31.1 & $10^{-5}$--$10^{10}$ & 0.0-27.3 & $10^{5}$--$10^{10}$& 0.0-47.7\\
     \bottomrule 
\end{tabular}
    \label{Pen_Table}
    \vspace{-0.65cm}
\end{table}
 We investigate the performance of the three penalization methods on 45 PGLib OPF test cases with up to 300 buses. Of these, the SDP relaxation is exact for $22.2$\%, i.e., an AC-feasible and globally optimal solution can be recovered. For the remaining $77.8$\% of test cases, we evaluate a wide range of penalty weights from $\epsilon_{\text{pen}} = \{ 10^{-5}, 10^{-4}, ..., 10^{9}, 10^{10} \}$\% and determine the number of additional test cases \bl{which lead to an AC-feasible solution}. For each successful test case, we evaluate the range of minimum and maximum penalty weights $\epsilon_{\text{pen}}^{\text{min}}$, $\epsilon_{\text{pen}}^{\text{max}}$ that allow recovery of an AC-feasible solution. For both the minimum and maximum penalty, we report the range of  minimum and maximum optimality gap with respect to the non-penalized objective value $f_{\text{cost}}^0$ of the SDP relaxation. 
 
 The results are shown in Table~\ref{Pen_Table}. Note that we display a range of values over all 35 investigated test cases (15 test cases are AC-feasible without penalization). The penalty term for reactive power is the most effective at recovering rank-1 solutions. Specifically, AC-feasible solutions are obtained for an additional $42.2$\% test cases. Note that the test cases recovered by the apparent branch loss and the matrix trace are a subset of those recovered by the reactive power penalty. As a result, in $35.6$\% of the test cases, none of these penalties successfully recover a rank-1 solution and the penalization heuristics fail to obtain an AC-feasible solution. \bl{For the reactive power penalty in particular, the spread of the minimum optimality gap over all test cases to obtain an AC-feasible solution is $24.0$\%, ranging from $0.0$\% to $24.0$\%. The sub-optimality for the reactive power penalty can increase up to $39.2\%$, indicating that a large sub-optimality can be incurred by assigning a sub-optimal penalty weight. Similar observations hold true for the apparent branch flow and matrix trace penalization terms.} Note the minimum penalty weight necessary to recover an AC-feasible solution for the different test cases and different penalty terms varies considerably, with the interval for the minimum reactive power penalty $\epsilon_{\text{pen}}^{\text{min}}$ ranging from $10^{-5}\%$ to $10^9\%$. Since a detailed screening of a wide range of penalty terms is likely to be computationally prohibitive, these results highlight the challenge of choosing a penalty weight that both recovers an AC-feasible solution and is small enough to obtain a near-globally optimal solution. Some works (e.g. \cite{Madani2016}) propose using a combination of penalty terms with individual penalty weights to obtain an AC-feasible solution, which, however, leads to an exponential increase in possible penalty weight combinations.
 
 \begin{figure}
    \centering
    \begin{footnotesize}
%
%
\definecolor{mycolor1}{rgb}{0.00000,0.44700,0.74100}%
\definecolor{mycolor2}{rgb}{0.85000,0.32500,0.09800}%
\begin{tikzpicture}

\begin{axis}[%
width=4.2cm,
height=5cm,
xmode=log,
xmin=0.001,
xmax=100,
xtick = {0.01, 0.1, 1,10,100},
xticklabels = { 0.01, 0.1, 1,10,100},
xlabel style={font=\color{white!15!black}},
xlabel={\textbf{a)} -- Optimality gap (\%)},
ymin=0.00001,
ymax=100,
ytick = {0.00001,0.0001,0.001, 0.01, 0.1, 1,10,100},
yticklabels={$10^{-5}$,$10^{-4}$,$10^{-3}$,$10^{-2}$,$10^{-1}$, $10^0$, $10^1$, $10^2$},
ymode = log,
ylabel style={font=\color{white!15!black}},
ylabel={Average distance to local optimality (\%)},
axis background/.style={fill=white},
xminorticks=true,
xmajorgrids = true,
ymajorgrids = true,
legend style={legend pos = north west, legend cell align=left, align=left, draw=white!15!black}
]
\addplot[only marks, mark=o, mark options={}, mark size=2pt, draw=mycolor1] table[row sep=crcr]{%
x	y\\
0	0.000411346839999601\\
0	0.00249387374732101\\
0	0.3029059160772\\
0	0.185337234405924\\
0	0.0292077685776204\\
0.00609675015144084	0.000500329655273987\\
0.199360286160943	2.66077183001708e-05\\
0	0.0125973192760021\\
0.0856412043132471	1.68120900267944\\
0.203877715753586	1.4769886789312\\
0	0.0726335394238293\\
0	0.040834313219751\\
24.0062928501776	0.0827996928999614\\
0.452775255343729	0.512714487642256\\
0.20576543084152	0.0198480978899487\\
0.395874232991411	0.692946651604535\\
0.0141338945124581	0.0287795953576292\\
0.044486425772261	0.00942631910072308\\
0.0380193464962542	5.48326363696302e-05\\
2.62242737407685	0.510425327091657\\
0.311460200994085	3.07625671115028\\
0.0292389071476995	0.294831190555163\\
0	0.0162975932093059\\
0.0781889065134989	1.55089651770517\\
0.127730394014214	0.339222250787749\\
2.22980185417858	2.89424839749705\\
0.110680472462721	1.44697545596393\\
8.20234431537701	9.94719345470065\\
0	0.245785046697575\\
};

\addplot[only marks, mark=square, mark options={}, mark size=2pt, draw=mycolor2] table[row sep=crcr]{%
x	y\\
0	0.000411346839999601\\
0	0.00249387374732101\\
0	0.3029059160772\\
0	0.185337234405924\\
0	0.0292077685776204\\
0.00617005872820542	0.00204160496078988\\
0.199360286160943	2.66077183001708e-05\\
0	0.0125973192760021\\
0.0871275051290255	1.70069403349654\\
0.192592096864574	0.606413473432869\\
0	0.0726335394238293\\
0	0.040834313219751\\
27.3370068924739	0.757220226452806\\
1.14658976930062	0.865979315160256\\
0.205757042607968	0.0165683096303835\\
0.0141099544005874	0.0128152646864225\\
0.044486425772261	0.00942631910072308\\
0.0380165602396287	0.000138073562356633\\
15.212617114458	7.4570568922624\\
0.209752065220914	1.45375347889627\\
0.0296761905012266	0.415752967059064\\
0	0.0162975932093059\\
0.106068812730742	1.38747507297298\\
0	0.245785046697575\\
};

\addplot[only marks, mark=diamond, mark options={}, mark size=2pt, draw=black] table[row sep=crcr]{%
x	y\\
0	0.000411346839999601\\
0	0.00249387374732101\\
0	0.3029059160772\\
0	0.185337234405924\\
0	0.0292077685776204\\
0.159081029843211	17.7799266534508\\
0.199360286160943	2.66077183001708e-05\\
0	0.0125973192760021\\
0	0.0726335394238293\\
0	0.040834313219751\\
24.0030675383678	0.0805990826438077\\
0.205760545565092	0.00118997522599512\\
0.0189659358319494	0.396376121542438\\
0.044486425772261	0.00942631910072308\\
0.0380186089573575	3.65511229313889e-05\\
0.553544500428416	8.0893506742354\\
0	0.0162975932093059\\
0	0.245785046697575\\
};

\end{axis}
\end{tikzpicture}%
 \, 
%
%
\definecolor{mycolor1}{rgb}{0.00000,0.44700,0.74100}%
\definecolor{mycolor2}{rgb}{0.85000,0.32500,0.09800}%
\begin{tikzpicture}

\begin{axis}[%
width=4.2cm,
height=5cm,
xmode=log,
xmin=1e-01,
xmax=1e04,
xtick = {1e-01, 1e00, 1e01, 1e02, 1e03,1e04},
xlabel style={font=\color{white!15!black}},
xlabel={\textbf{b)} -- Distance to AC-feasiblity (\%)},
ymin=0.01,
ymax=100,
ytick = {0.01,0.1, 1, 10, 100},
yticklabels={$10^{-2}$,$10^{-1}$, $10^{0}$, $10^{1}$, $10^{2}$},
ymode=log,
axis background/.style={fill=white},
xminorticks=true,
xmajorgrids = true,
ymajorgrids = true,
scaled y ticks=false,
legend style={legend pos = north west, legend cell align=left, align=left, draw=white!15!black}
]
\addplot[only marks, mark=o, mark options={}, mark size=2pt, draw=mycolor1] table[row sep=crcr]{%
x	y\\
181.710739681513	6.50173511870853\\
549.079061698111	7.65842427705574\\
108.187157126049	12.0841085441985\\
4.50624551637729	6.63261305633471\\
0.624146232526767	9.72118332511268\\
1.72837394404615	11.6644150530329\\
4.30806382986905	7.20253966794652\\
5.6066147792181	12.31306197422\\
151.356828233186	3.76098697277427\\
598.146457710445	11.3696573443086\\
113.915131097034	8.49630374528597\\
4.04526153135409	8.66292736088886\\
223.087437593505	6.81483001132334\\
1028.08055193116	5.97961954573823\\
181.823827443017	9.10723603018834\\
4.47343877093556	6.46587603963162\\
};

\addplot[only marks, mark=square, mark options={}, mark size=2pt, draw=mycolor2] table[row sep=crcr]{%
x	y\\
745.302511465816	6.39270383685479\\
151.217150925613	5.40536079356941\\
342.00451949648	9.62945869774389\\
16.6603101747669	0.582300042406268\\
8.34685380641842	1.92138763926827\\
3.92319148998386	1.60997765800734\\
15.2964123831896	4.99618295752587\\
4.60385476445378	7.07267992312801\\
8.30256262239014	10.6741364818109\\
249.010306037229	5.25816452830298\\
832.653586140327	8.03891281349449\\
284.403217806472	6.50397037585372\\
10.7628484585372	1.02062717285919\\
1.33327917723631	0.173447277388062\\
0.281524056919107	0.0910202085555016\\
4.97385144304372	6.63146059242204\\
2.11620152059185	12.5756591042417\\
1005.56656096402	6.85957009806498\\
1651.31173768204	9.53790631851289\\
441.509511035295	9.43923861792124\\
11.8405387718984	3.6543037027291\\
};

\addplot[only marks, mark=diamond, mark options={}, mark size=2pt, draw=black] table[row sep=crcr]{%
x	y\\
0.383931164783929	5.51795919679995\\
0.138606177946717	21.3649088395505\\
178.290196575727	7.10801034756287\\
7.34672135991828	14.0833379962015\\
344.131653473493	10.6608470123401\\
10.5290215109273	10.3283192516422\\
19.0039603462524	12.1390584137188\\
1.17259685936362	0.906998212086156\\
3.92319148998386	1.60997765800734\\
13.9125480211804	7.98644407526859\\
7.47163399246979	13.0155839331288\\
40.8950253180751	14.4226729802985\\
41.8166196806806	3.62063347080653\\
376.848469385289	8.31771310585159\\
234.531091616033	12.457645259025\\
4.32710025444344	14.0752398665433\\
5.19254346446697	18.0201096443526\\
1.87382859424272	8.47675030435169\\
1.50993215498217	6.31072630615904\\
1.50450849305176	0.513993887489893\\
18.0651091275055	18.6806607376764\\
0.445811664388389	2.39140708978879\\
12.4014980555768	19.7267408166683\\
144.657822208135	4.99632846324336\\
66.9797246414098	14.761950932249\\
436.796782925136	9.44247997527315\\
11.3118073101646	10.380068475613\\
};

\end{axis}
\end{tikzpicture}%

    \definecolor{mycolor1}{rgb}{0.00000,0.44700,0.74100}%
\definecolor{mycolor2}{rgb}{0.85000,0.32500,0.09800}%
    \begin{tikzpicture} 
    \begin{axis}[%
    hide axis,
    xmin=1,
    xmax=2,
    ymin=0,
    ymax=0.1,
    legend style={draw=white!15!black,legend cell align=left},
    legend columns=3
    ]
    \addlegendimage{only marks, mark=o, mark options={}, mark size=3pt, draw=mycolor1}
    \addlegendentry{Reactive power $Q_G$};
            \addlegendimage{only marks, mark=diamond, mark options={}, mark size=3pt, draw=black}
    \addlegendentry{Matrix trace $\text{Tr}\{W\}$};
        \addlegendimage{only marks, mark=square, mark options={}, mark size=3pt, draw=mycolor2}
    \addlegendentry{Branch flow $| S_{ij} - S_{ji}|$};
    \end{axis}
\end{tikzpicture}
\end{footnotesize}
   \vspace{-0.75cm}
    \caption{\bl{We show a) the optimality gap versus the average normalized distance to local optimality for test cases where the penalization has been successful, and b) the cumulative normalized constraint violation versus the average normalized distance to a locally optimal solution for test cases where the penalization has failed. Both axes are logarithmic. Note that for a) we report the metrics for the minimum (non-zero) magnitude of penalty term $\epsilon_{\text{min}}$ necessary to achieve an AC-feasible solution. For b) we report the metrics for the minimum cumulative normalized constraint violation, i.e., the magnitude of penalty term for which the solution variables are closest to AC-feasibility.}}
    \label{Pen_Detail}
    \vspace{-0.45cm}
\end{figure}
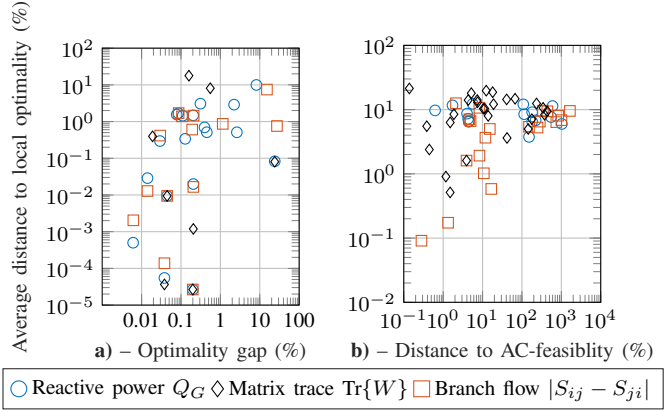
 \bl{To provide more insight into the quality of the solutions obtained from the three different penalization techniques, we show in Fig.~\ref{Pen_Detail} a) the optimality gap versus the average normalized distance to local optimality for test cases where the penalization has been successful, and b) the cumulative normalized constraint violation versus the average normalized distance to a locally optimal solution for test cases where the penalization failed. Note that for a) we report the metrics for the minimum (non-zero) magnitude of penalty term $\epsilon_{\text{min}}$ necessary to achieve an AC-feasible solution. For b) we report the metrics for the minimum value of the distance to AC-feasibility in \eqref{eq:constrviol}, i.e., the magnitude of penalty term for which the solution variables are closest to AC-feasibility.}
 
  \bl{First, focusing on Fig.~\ref{Pen_Detail}a), we observe the wide spread between the minimum optimality gaps for the different test cases and penalization terms (corresponding to the third column in Table~\ref{Pen_Table}). We observe that the average distance to the locally optimal solution provided by a non-convex solver can also be substantial (larger than 10\%). This confirms the findings of the 5-bus test case in the previous subsection that the penalization terms do not necessarily drive the solution towards the locally optimal solution, but instead toward other parts of the non-convex feasible space with higher costs. This highlights that AC-feasible solutions obtained by penalized SDP relaxations might incur substantial sub-optimality compared to the solution obtained by a non-convex solver.} 
 
  Second, focusing on Fig.~\ref{Pen_Detail}b), we can observe that for the majority of cases for which at least one of the penalization methods did not achieve an AC-feasible solution, both the distance to AC-feasibility and distance to the locally optimal solution is substantial. This showcases that the penalty terms can also drive the solution towards portions of the feasible space of the convex relaxation that are not close to the feasible space of the non-convex AC-OPF. Only 5 out of the 64 combinations of test cases and penalty terms have a distance to AC-feasibility that is below 1\%. For these 5 test cases only, we re-run the penalization with penalization weights sampled close to the penalty weight that lead to the smallest cumulative normalized constraint violation in \eqref{eq:constrviol}. As a result, we identify penalty weights for these test cases that lead to an AC-feasible solution. To conclude, applying the proposed metrics for penalization methods allows us to fine-tune the penalty weights to obtain AC-feasible solutions. On the other hand, as evident in Fig.~\ref{Pen_Detail}b), there exist several test cases (33.3\% of test cases examined) for which none of the penalization methods resulted in recovery of an AC-feasible solution. The penalized solutions have substantial distances both to AC-feasibility and to the locally optimal solution provided by a non-convex solver. This motivates the investigation of warm-starting non-convex solvers with the inexact solutions of convex relaxations to recover an AC-feasible solution in the next section. \blueEPSR{We compute both the Spearman's rank and Pearson's correlation coefficient to quantify the correlations in Fig.~\ref{Pen_Detail}. The correlation coefficients range from -0.28 to 0.53, thus showing that neither of the quantities have a strong correlation using both correlation coefficients.}
 
\subsection{Warm-Starting Non-Convex Solvers}
\label{Warm_Start}
This section investigates whether non-convex solvers can be efficiently warm-started in the search for local optima when initialized with solutions of convex relaxations. To this end, we use the SQP solver in KNITRO (algorithm 4), and the IPM solvers provided by KNITRO (algorithm 1) and IPOPT. We deactivated the presolve in KNITRO, as enabling the presolve resulted in significantly longer solver times. An upper time limit of 2000 seconds is enforced. For the remaining options, we use the default values. We first look at the solver reliability, then study the variation in computational speed for different initializations, and finally evaluate the solution quality. 
 \subsubsection{Solver Reliability}
\begin{table}[]
\caption{PGLib OPF test cases solved to local optimality (\%)}
    \centering
\begin{tabular}{c c c c c}
\toprule 
 
   \quad  Solver / Initialization & Flat start & DC-OPF & QC & SDP \\
     \midrule 
         KNITRO (SQP) & 85.4 & 75.0 & 81.3 & 81.3 \\
        KNITRO (IPM) & 99.0 & 92.7 & 97.9 & 93.8 \\
          IPOPT (IPM) & 100.0 & 100.0 & 100.0 & 100.0 \\

     \bottomrule
\end{tabular}
    \label{Feas_Table}
    \vspace{-0.4cm}
\end{table}
Table~\ref{Feas_Table} shows the share of the 96 considered PGLib OPF test cases which are solved to local optimality for the different initializations and solvers. IPOPT is the most reliable solver, with 100\% of the test cases solved to local optimality irrespective of the initialization. The IPM and SQP solvers in KNITRO are less reliable, and achieve their highest reliability for the flat start initialization. For the other instances solved by KNITRO, either the time limit of 2000 seconds was reached or the solver reported local infeasibility. 
 \subsubsection{Computational Speed}
 \begin{figure}
    \centering
    \begin{footnotesize}
%
%
\begin{tikzpicture}

\begin{axis}[%
width=3.3cm,
height=4.3cm,
xmin=0.5,
xmax=3.5,
xtick={1,2,3},
xticklabels={{DC},{QC},{SDP}},
ymin=3,
ymax=1000,
ymode=log,
ylabel style={font=\color{white!15!black}},
ylabel={Speed -- KNITRO, SQP (\%)},
axis background/.style={fill=white},
xmajorgrids,
ymajorgrids
]
\addplot [color=black, dotted, line width=1pt,forget plot]
  table[row sep=crcr]{%
-1	100\\
4	100\\
};
\addplot [color=black, dashed, forget plot]
  table[row sep=crcr]{%
1	98.1743805721166\\
1	131.57936761712\\
};
\addplot [color=black, dashed, forget plot]
  table[row sep=crcr]{%
2	104.252538660554\\
2	150.263780643803\\
};
\addplot [color=black, dashed, forget plot]
  table[row sep=crcr]{%
3	112.60626512395\\
3	141.035533389276\\
};
\addplot [color=black, dashed, forget plot]
  table[row sep=crcr]{%
1	51.5067973460545\\
1	72.8636149627628\\
};
\addplot [color=black, dashed, forget plot]
  table[row sep=crcr]{%
2	38.3022624094886\\
2	72.84957606805\\
};
\addplot [color=black, dashed, forget plot]
  table[row sep=crcr]{%
3	52.5488015999044\\
3	80.8812691321017\\
};
\addplot [color=black, forget plot]
  table[row sep=crcr]{%
0.8875	131.57936761712\\
1.1125	131.57936761712\\
};
\addplot [color=black, forget plot]
  table[row sep=crcr]{%
1.8875	150.263780643803\\
2.1125	150.263780643803\\
};
\addplot [color=black, forget plot]
  table[row sep=crcr]{%
2.8875	141.035533389276\\
3.1125	141.035533389276\\
};
\addplot [color=black, forget plot]
  table[row sep=crcr]{%
0.8875	51.5067973460545\\
1.1125	51.5067973460545\\
};
\addplot [color=black, forget plot]
  table[row sep=crcr]{%
1.8875	38.3022624094886\\
2.1125	38.3022624094886\\
};
\addplot [color=black, forget plot]
  table[row sep=crcr]{%
2.8875	52.5488015999044\\
3.1125	52.5488015999044\\
};
\addplot [color=blue, forget plot]
  table[row sep=crcr]{%
0.775	72.8636149627628\\
0.775	98.1743805721166\\
1.225	98.1743805721166\\
1.225	72.8636149627628\\
0.775	72.8636149627628\\
};
\addplot [color=blue, forget plot]
  table[row sep=crcr]{%
1.775	72.84957606805\\
1.775	104.252538660554\\
2.225	104.252538660554\\
2.225	72.84957606805\\
1.775	72.84957606805\\
};
\addplot [color=blue, forget plot]
  table[row sep=crcr]{%
2.775	80.8812691321017\\
2.775	112.60626512395\\
3.225	112.60626512395\\
3.225	80.8812691321017\\
2.775	80.8812691321017\\
};
\addplot [color=red, forget plot]
  table[row sep=crcr]{%
0.775	83.37793264756\\
1.225	83.37793264756\\
};
\addplot [color=red, forget plot]
  table[row sep=crcr]{%
1.775	89.0756411956831\\
2.225	89.0756411956831\\
};
\addplot [color=red, forget plot]
  table[row sep=crcr]{%
2.775	94.5055434502574\\
3.225	94.5055434502574\\
};
\addplot [color=black, draw=none, mark=+, mark options={solid, red}, forget plot]
  table[row sep=crcr]{%
1	144.953170528494\\
1	339.587576719032\\
};
\addplot [color=black, draw=none, mark=+, mark options={solid, red}, forget plot]
  table[row sep=crcr]{%
2	22.9800393950993\\
2	162.748855420123\\
2	208.080225947797\\
2	235.400267749523\\
2	271.928738958116\\
2	332.443536462921\\
2	461.255233075591\\
};
\addplot [color=black, draw=none, mark=+, mark options={solid, red}, forget plot]
  table[row sep=crcr]{%
3	26.3051161405475\\
3	162.541999029462\\
3	171.897645624609\\
3	430.332872692287\\
3	983.540448248642\\
};
\end{axis}
\end{tikzpicture}%
%
%
\begin{tikzpicture}

\begin{axis}[%
width=3.3cm,
height=4.3cm,
xmin=0.5,
xmax=3.5,
xtick={1,2,3},
xticklabels={{DC},{QC},{SDP}},
ymin=3,
ymax=1000,
ymode=log,
ylabel style={font=\color{white!15!black}},
ylabel={Speed -- KNITRO, IPM (\%)},
axis background/.style={fill=white},
xmajorgrids,
ymajorgrids
]
\addplot [color=black, dotted, line width=1pt,forget plot]
  table[row sep=crcr]{%
-1	100\\
4	100\\
};
\addplot [color=black, dashed, forget plot]
  table[row sep=crcr]{%
1	112.663792405892\\
1	155.555849901321\\
};
\addplot [color=black, dashed, forget plot]
  table[row sep=crcr]{%
2	147.391157089753\\
2	230.608368207393\\
};
\addplot [color=black, dashed, forget plot]
  table[row sep=crcr]{%
3	140.001636337147\\
3	214.498512224575\\
};
\addplot [color=black, dashed, forget plot]
  table[row sep=crcr]{%
1	43.4210567593258\\
1	82.5622164121184\\
};
\addplot [color=black, dashed, forget plot]
  table[row sep=crcr]{%
2	4.37497857959227\\
2	77.1859870498685\\
};
\addplot [color=black, dashed, forget plot]
  table[row sep=crcr]{%
3	4.9999329448106\\
3	74.3172571328541\\
};
\addplot [color=black, forget plot]
  table[row sep=crcr]{%
0.8875	155.555849901321\\
1.1125	155.555849901321\\
};
\addplot [color=black, forget plot]
  table[row sep=crcr]{%
1.8875	230.608368207393\\
2.1125	230.608368207393\\
};
\addplot [color=black, forget plot]
  table[row sep=crcr]{%
2.8875	214.498512224575\\
3.1125	214.498512224575\\
};
\addplot [color=black, forget plot]
  table[row sep=crcr]{%
0.8875	43.4210567593258\\
1.1125	43.4210567593258\\
};
\addplot [color=black, forget plot]
  table[row sep=crcr]{%
1.8875	4.37497857959227\\
2.1125	4.37497857959227\\
};
\addplot [color=black, forget plot]
  table[row sep=crcr]{%
2.8875	4.9999329448106\\
3.1125	4.9999329448106\\
};
\addplot [color=blue, forget plot]
  table[row sep=crcr]{%
0.775	82.5622164121184\\
0.775	112.663792405892\\
1.225	112.663792405892\\
1.225	82.5622164121184\\
0.775	82.5622164121184\\
};
\addplot [color=blue, forget plot]
  table[row sep=crcr]{%
1.775	77.1859870498685\\
1.775	147.391157089753\\
2.225	147.391157089753\\
2.225	77.1859870498685\\
1.775	77.1859870498685\\
};
\addplot [color=blue, forget plot]
  table[row sep=crcr]{%
2.775	74.3172571328541\\
2.775	140.001636337147\\
3.225	140.001636337147\\
3.225	74.3172571328541\\
2.775	74.3172571328541\\
};
\addplot [color=red, forget plot]
  table[row sep=crcr]{%
0.775	95.3580691471596\\
1.225	95.3580691471596\\
};
\addplot [color=red, forget plot]
  table[row sep=crcr]{%
1.775	99.9995231646608\\
2.225	99.9995231646608\\
};
\addplot [color=red, forget plot]
  table[row sep=crcr]{%
2.775	102.157495139039\\
3.225	102.157495139039\\
};
\addplot [color=black, draw=none, mark=+, mark options={solid, red}, forget plot]
  table[row sep=crcr]{%
1	4.37497857959227\\
1	18.1318479781598\\
1	20.7302038581115\\
1	26.7748382065434\\
1	172.221780709423\\
1	368.354640253757\\
};
\addplot [color=black, draw=none, mark=+, mark options={solid, red}, forget plot]
  table[row sep=crcr]{%
2	256.043636926456\\
2	260.616286343421\\
2	268.391470771455\\
2	269.126968610301\\
2	271.2706191189\\
2	274.852028850054\\
2	296.654492131811\\
2	301.426033184382\\
2	347.146068613886\\
2	650.901493340153\\
2	663.22247151689\\
2	758.142360353591\\
2	968.826307418285\\
};
\addplot [color=black, draw=none, mark=+, mark options={solid, red}, forget plot]
  table[row sep=crcr]{%
3	275.003476982377\\
3	297.757010774669\\
3	301.32470827846\\
3	310.067811530392\\
3	362.52317858287\\
3	1484.5964124711\\
};
\end{axis}
\end{tikzpicture}%
%
%
\begin{tikzpicture}

\begin{axis}[%
width=3.3cm,
height=4.3cm,
xmin=0.5,
xmax=3.5,
ymode=log,
xtick={1,2,3},
xticklabels={{DC},{QC},{SDP}},
ymin=3,
ymax=1000,
ylabel style={font=\color{white!15!black}},
ylabel={Speed -- IPOPT, IPM (\%)},
axis background/.style={fill=white},
xmajorgrids,
ymajorgrids
]
\addplot [color=black, dotted, line width=1pt,forget plot]
  table[row sep=crcr]{%
-1	100\\
4	100\\
};
\addplot [color=black, dashed, forget plot]
  table[row sep=crcr]{%
1	87.7623015760671\\
1	118.880961810283\\
};
\addplot [color=black, dashed, forget plot]
  table[row sep=crcr]{%
2	100.972535492444\\
2	173.823269669713\\
};
\addplot [color=black, dashed, forget plot]
  table[row sep=crcr]{%
3	101.725495732417\\
3	159.534292788511\\
};
\addplot [color=black, dashed, forget plot]
  table[row sep=crcr]{%
1	17.5200342918595\\
1	59.3431934184667\\
};
\addplot [color=black, dashed, forget plot]
  table[row sep=crcr]{%
2	5.19093058608068\\
2	51.9800343314517\\
};
\addplot [color=black, dashed, forget plot]
  table[row sep=crcr]{%
3	5.36759971492846\\
3	59.9791082373345\\
};
\addplot [color=black, forget plot]
  table[row sep=crcr]{%
0.8875	118.880961810283\\
1.1125	118.880961810283\\
};
\addplot [color=black, forget plot]
  table[row sep=crcr]{%
1.8875	173.823269669713\\
2.1125	173.823269669713\\
};
\addplot [color=black, forget plot]
  table[row sep=crcr]{%
2.8875	159.534292788511\\
3.1125	159.534292788511\\
};
\addplot [color=black, forget plot]
  table[row sep=crcr]{%
0.8875	17.5200342918595\\
1.1125	17.5200342918595\\
};
\addplot [color=black, forget plot]
  table[row sep=crcr]{%
1.8875	5.19093058608068\\
2.1125	5.19093058608068\\
};
\addplot [color=black, forget plot]
  table[row sep=crcr]{%
2.8875	5.36759971492846\\
3.1125	5.36759971492846\\
};
\addplot [color=blue, forget plot]
  table[row sep=crcr]{%
0.775	59.3431934184667\\
0.775	87.7623015760671\\
1.225	87.7623015760671\\
1.225	59.3431934184667\\
0.775	59.3431934184667\\
};
\addplot [color=blue, forget plot]
  table[row sep=crcr]{%
1.775	51.9800343314517\\
1.775	100.972535492444\\
2.225	100.972535492444\\
2.225	51.9800343314517\\
1.775	51.9800343314517\\
};
\addplot [color=blue, forget plot]
  table[row sep=crcr]{%
2.775	59.9791082373345\\
2.775	101.725495732417\\
3.225	101.725495732417\\
3.225	59.9791082373345\\
2.775	59.9791082373345\\
};
\addplot [color=red, forget plot]
  table[row sep=crcr]{%
0.775	74.1525571899635\\
1.225	74.1525571899635\\
};
\addplot [color=red, forget plot]
  table[row sep=crcr]{%
1.775	76.1057072995986\\
2.225	76.1057072995986\\
};
\addplot [color=red, forget plot]
  table[row sep=crcr]{%
2.775	82.6022038487846\\
3.225	82.6022038487846\\
};
\addplot [color=black, draw=none, mark=+, mark options={solid, red}, forget plot]
  table[row sep=crcr]{%
1	5.02286939379543\\
1	14.3974764333259\\
1	137.704144070688\\
1	146.89085821138\\
1	193.344123998061\\
1	210.879296968598\\
1	221.857678743145\\
1	253.997435695083\\
};
\addplot [color=black, draw=none, mark=+, mark options={solid, red}, forget plot]
  table[row sep=crcr]{%
2	193.885976251332\\
2	194.062650935127\\
2	197.373761710518\\
2	200.630204667662\\
2	203.81520187552\\
2	240.746219449771\\
2	260.335515004449\\
2	267.085245560809\\
2	331.083299618397\\
2	360.329845797916\\
2	435.741649779041\\
2	452.95001759487\\
2	943.375087698815\\
};
\addplot [color=black, draw=none, mark=+, mark options={solid, red}, forget plot]
  table[row sep=crcr]{%
3	209.848243546064\\
3	231.384983397978\\
3	245.130896660395\\
3	339.851447559757\\
3	424.955451698785\\
};
\end{axis}
\end{tikzpicture}%
    \end{footnotesize}
    \vspace{-0.65cm}
    \caption{Variation in computational speed relative to a flat start resulting from warm-starting the SQP solver in KNITRO, the IPM solver in KNITRO, and IPOPT using the solutions of the DC-OPF, the QC relaxation, and the SDP relaxation. Note the y-axis is shown on a logarithmic scale. The value of $100$\% ($10^2$\%) corresponds to the computational speed of the flat start, with lower values indicating a speed improvement and higher values slower performance.}
    \label{IP_Speed}
    \vspace{-0.4cm}
\end{figure}
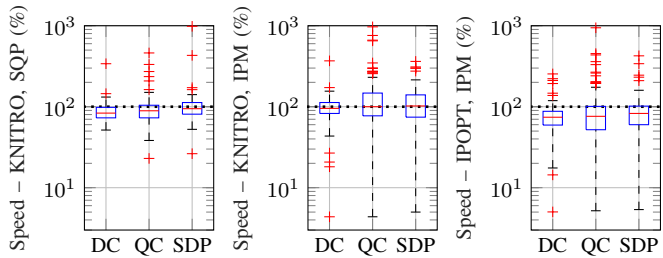
Fig.~\ref{IP_Speed} shows the variation in computational speed relative to a flat start for the warm-started non-convex solvers initialized with the solutions to the DC-OPF and the QC and SDP relaxations. Note that we only consider instances solved to local optimality, and we do not include the computational time required to compute the initializations. Warm-starting can have a positive or negative effect on computational speed for both solution methods and all three solvers. \bl{The interior-point solver IPOPT shows the best performance with the lower 75th percentiles exhibiting speed improvements when initalized with solutions to the DC-OPF and the QC and SDP relaxations as well as median speed improvements of $25.8$\% for DC-OPF, $23.9$\% for the QC relaxation, and $17.4$\% for the SDP relaxation. This does not confirm that SQP methods are usually more suitable for warm-starting as stated in Section~\ref{Warm_start_theory}.} The IPM solver in KNITRO performs worse, with only the lower 50th percentile of test cases exhibiting a speed improvement. The SQP solver in KNITRO has better computational speed than the IPM solver in KNITRO but has significantly lower solver reliability as shown in Table~\ref{Feas_Table}. \blue{Both the Pearson's and Spearman's rank correlation coefficients for the computational speed-up and the i) optimality gaps, ii) cumulative constraint violations, and iii) the distances to local optimality lie in a range between -0.22 and 0.49, thus showing that these three metrics are not strongly correlated with the computational speed-up in these statistical measures.}

\begin{figure}
    \centering
    \begin{footnotesize}
    \input{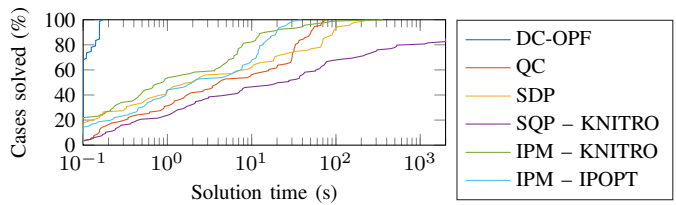}
    \end{footnotesize}
    \vspace{-0.75cm}
    \caption{\bl{PGLib OPF instances solved versus computational time for the DC-OPF, the QC and SDP relaxations as well as for the SQP solver in KNITRO, the IPM solver in KNITRO, and IPOPT when initialized with a flat start.}}
    \label{IPOPT_vs_SQP}
    \vspace{-0.4cm}
\end{figure}
\bl{Fig.~\ref{IPOPT_vs_SQP} shows the solution times of all the non-convex solvers initalized with a flat start and the solution times for the DC-OPF, the QC relaxation, and the SDP relaxation. The DC-OPF is by far the fastest, as it involves only solving an LP or QP. The SDP relaxation in PowerModels.jl is faster than the QC relaxation for small and medium systems but is approximately an order-of-magnitude slower for large systems. It can be observed that interior-point solvers, in particular IPOPT, are significantly faster than SQP solvers for the AC-OPF problem as the SQP solvers do not scale well with increasing system size.} \blueEPSR{Note that the computational time for solving the QC and SDP relaxation (as shown in Fig.~\ref{IPOPT_vs_SQP}) is larger than the computational time for solving the AC-OPF with an interior-point solver and a flat start initialization in most instances. As a result, even if a very large speed-up is achieved by initializing the non-convex solver with the solution to the QC and SDP relaxation, the overall combined computational time still will be larger than directly solving the AC-OPF with a flat start. This finding highlights that there is a need for computationally efficient methods to recover an AC-feasible solution from inexact convex relaxations.}

\subsubsection{Solution Quality}
For the evaluated cases, all non-convex solvers which converge to local optimality from all starting points obtain the same objective value to within a small numerical tolerance of $ 10^{-5}$. \rev{This confirms the finding of the study \cite{kardos2018complete} that non-convex solvers report the same objective value for a wide range of cases and different spurious locally optimal solutions are not identified.} For the five-bus test case from \cite{Bukhsh2013} discussed earlier, IPOPT and the SQP solver in KNITRO return the spurious locally optimal solution, i.e. not the globally optimal one, when initialized with the flat start. All warm-started local solvers return the globally optimal solution to this problem; \blue{thus, in this case, warm-starting the solvers yields an improvement in solution quality.}t

\section{Conclusions and Outlook}
\label{sec:conclusion}
Using the PGLib OPF benchmarks from \cite{PGLIB}, we provided a comprehensive study regarding the distance of the inexact solutions to the QC relaxation and the SDP relaxation relative to both AC feasibility and local optimality for the original non-convex AC-OPF problem. We investigated penalization methods for recovering AC-feasible solutions and warm-starting of non-convex solvers for recovering locally optimal solutions. Based on our detailed results, we summarize our main conclusions and outline directions for further research:

1)~To quantify the distances to AC feasibility and local optimality, we proposed two empirical metrics to complement the optimality gap for analyzing the relaxations' accuracy. For both the QC and SDP relaxations, we have shown that these metrics are not strongly correlated with the optimality gap, and despite an optimality gap of less than $1$\%, several test cases still exhibit a substantial distance to AC feasibility and local optimality \rev{highlighting the added value of the two metrics.} Detailed investigations of these cases and of four outliers could provide additional insights into the QC and SDP relaxations. 
    
    2)~Heuristic penalization methods for the SDP relaxation can be successful in recovering an AC-feasible and near-globally optimal solution. However, choosing an appropriate penalty weight can be challenging since detailed screenings of a wide range of penalty terms might be computationally prohibitive. Furthermore, there exists a range of test cases for which all three penalization methods fail to recover an AC-feasible solution. \bl{We characterize the obtained solutions from the penalized SDP relaxation using our proposed metrics, and show that most instances for which penalization methods fail exhibit substantial distances to both AC-feasibility and local optimality. For failed test instances with small distances, we show how our proposed metrics inform a fine-tuning of penalty weights to obtain AC-feasible solutions.} \blueEPSR{A direction for improvement is to develop systematic and scalable methods to choose penalty weights and penalty terms in order to reliably recover AC-feasible and near-globally optimal solutions.} 
    
            3)~We investigated warm-starting non-convex solvers with the solutions of the inexact convex relaxations. We have shown that the sequential quadratic programming solver in KNITRO is not scalable to large instances and exhibits issues with solver convergence when warm started. Conversely, the interior-point solver in IPOPT is highly reliable and computationally efficient: in more than 75\% of the considered PGLib OPF test cases, a computational speed-up can be achieved by warm-starting with the solution of an inexact convex relaxation. A future direction is to improve the warm-starting of the interior-point method by a)~decreasing the initial logarithmic barrier term in the objective function, b)~using dual information or c) learning optimal warm-starts points \cite{baker2019learning}.
            \vspace{-0.325cm}
\bibliographystyle{IEEEtran}

\bibliography{Bib}

\appendix 
\blueEPSR{This appendix provides detailed results for the obtained optimality gaps as well as distances to AC-feasbility and local optimality for PGLib OPF test cases from \cite{PGLIB} in Table~\ref{Appendix_1}, Table~\ref{Appendix_2}, and Table~\ref{Appendix_3} for the base case benchmarks, heavily loaded test cases (api), and small phase angle difference cases (sad), respectively.}

\begin{table*}[]
    \centering
        \caption{Optimality gaps, distances to AC-feasibility and local optimality for PGLib OPF test cases from \cite{PGLIB}  (Base case benchmarks)}
    \label{Appendix_1}
\begin{tabular}{r c c c c c c}
\toprule
 & \multicolumn{3}{c}{QC relaxation} & \multicolumn{3}{c}{SDP relaxation} \\
\midrule
&  & distance to & distance to &  & distance to & distance to\\
 & opt. gap & AC-feas. & local opt.&  opt. gap  & AC-feas. & local opt. \\
Test cases & (\%) & (\%) & (\%) & (\%) &(\%) &(\%) \\
\midrule
 pglib\_opf\_case14\_ieee &  0.11 & 1.01e+00  &  1.30 &   0.00 &  0.00e+00 &   1.17\\ \midrule
          pglib\_opf\_case24\_ieee\_rts &  0.01 & 4.09e-01  &  0.71 &   0.00 &  4.11e-07 &   1.52\\ \midrule
                pglib\_opf\_case30\_as &  0.06 & 0.00e+00  &  1.34 &   0.00 &  0.00e+00 &   0.91\\ \midrule
               pglib\_opf\_case30\_fsr &  0.39 & 4.70e+00  &  8.11 &   0.00 &  8.62e-03 &   0.49\\ \midrule
              pglib\_opf\_case30\_ieee & 10.78 & 3.08e+01  &  3.35 &   0.00 &  9.41e-05 &   0.72\\ \midrule
              pglib\_opf\_case39\_epri &  0.48 & 3.47e+01  &  5.31 &   0.01 &  1.42e+01 &   1.43\\ \midrule
              pglib\_opf\_case57\_ieee &  0.46 & 3.08e+01  &  2.36 &   0.00 &  7.61e-06 &   1.14\\ \midrule
          pglib\_opf\_case73\_ieee\_rts &  0.03 & 3.71e+01  &  1.19 &   0.00 &  2.71e-05 &   1.53\\ \midrule
            pglib\_opf\_case89\_pegase &  0.73 & 3.26e+01  &  5.67 &   0.01 &  5.29e+01 &   1.71\\ \midrule
             pglib\_opf\_case118\_ieee &  2.19 & 2.62e+02  &  5.82 &   0.18 &  1.33e+02 &   4.09\\ \midrule
         pglib\_opf\_case162\_ieee\_dtc &  7.54 & 2.07e+02  &  7.68 &   2.26 &  1.02e+03 &   6.08\\ \midrule
              pglib\_opf\_case179\_goc &  0.12 & 2.04e+03  &  8.26 &   0.06 &  2.04e+03 &   7.98\\ \midrule
             pglib\_opf\_case200\_tamu &  0.00 & 5.27e-01  &  0.42 &   0.00 &  2.27e-02 &   0.29\\ \midrule
            pglib\_opf\_case240\_pserc &  3.80 & 7.30e+02  & 10.88 &   2.27 &  6.69e+02 &   9.80\\ \midrule
             pglib\_opf\_case300\_ieee &  2.55 & 2.13e+02  &  3.10 &   0.11 &  2.00e+01 &   1.78\\ \midrule
             pglib\_opf\_case500\_tamu &  5.38 & 7.43e+02  &  1.96 &   2.11 &  8.81e+03 &   3.39\\ \midrule
             pglib\_opf\_case588\_sdet &  1.67 & 8.20e+03  &  7.42 &   0.67 &  3.06e+02 &   4.37\\ \midrule
          pglib\_opf\_case1354\_pegase &  2.39 & 3.79e+03  &  3.96 &   0.56 &  1.62e+03 &   2.31\\ \midrule
             pglib\_opf\_case1888\_rte &  1.81 & n.a.  & 10.70 &   1.74 &  n.a. &  10.47\\ \midrule
             pglib\_opf\_case1951\_rte &  0.11 & n.a.  &  2.53 &   0.01 &  n.a. &   1.73\\ \midrule
            pglib\_opf\_case2316\_sdet &  2.21 & 7.81e+03  &  6.37 &   0.73 &  2.32e+03 &   4.51\\ \midrule
             pglib\_opf\_case2383wp\_k &  0.99 & 1.11e+04  &  4.60 &   0.38 &  5.12e+03 &   3.58\\ \midrule
             pglib\_opf\_case2736sp\_k &  0.29 & 1.22e+04  &  2.56 &   0.02 &  1.36e+03 &   0.36\\ \midrule
            pglib\_opf\_case2737sop\_k &  0.25 & 1.53e+04  &  1.87 &   0.03 &  8.48e+02 &   0.46\\ \midrule
            pglib\_opf\_case2746wop\_k &  0.36 & 2.99e+05  &  2.88 &   0.06 &  1.39e+04 &   0.52\\ \midrule
             pglib\_opf\_case2746wp\_k &  0.32 & 3.44e+04  &  2.55 &   0.01 &  1.71e+03 &   0.41\\ \midrule
             pglib\_opf\_case2848\_rte &  0.12 & n.a.  &  4.42 &   0.04 &  n.a. &   3.94\\ \midrule
            pglib\_opf\_case2853\_sdet &  1.43 & 5.30e+03  &  7.64 &   0.54 &  4.81e+03 &   5.64\\ \midrule
             pglib\_opf\_case2868\_rte &  0.10 & n.a.  &  2.16 &   0.20 &  n.a. &   2.29\\ \midrule
          pglib\_opf\_case2869\_pegase &  1.07 & 9.80e+03  &  2.86 &   0.41 &  1.28e+03 &   1.78\\ \midrule
             pglib\_opf\_case3012wp\_k &  0.98 & 5.24e+03  &  4.03 &   0.16 &  1.05e+04 &   0.99\\ \midrule
             pglib\_opf\_case3120sp\_k &  0.53 & 1.01e+04  &  3.35 &   0.11 &  5.20e+03 &   0.75\\ 
\bottomrule
\multicolumn{7}{l}{n.a. -- The AC power flows in MATPOWER did not converge for these test cases.} \\
\end{tabular}
\end{table*}

\begin{table*}[]
    \centering
        \caption{Optimality gaps, distances to AC-feasibility  and local optimality for PGLib OPF test cases from \cite{PGLIB}  (heavily loaded test cases (api))}
    \label{Appendix_2}
\begin{tabular}{r c c c c c c}
\toprule
 & \multicolumn{3}{c}{QC relaxation} & \multicolumn{3}{c}{SDP relaxation} \\
\midrule
&  & distance to & distance to &  & distance to & distance to\\
 & opt. gap & AC-feas. & local opt.&  opt. gap  & AC-feas. & local opt. \\
Test cases & (\%) & (\%) & (\%) & (\%) &(\%) &(\%) \\
\midrule
 pglib\_opf\_case14\_ieee\_\_api &  1.77 & 3.21e+00  &  1.13 &   0.00 &  1.62e-06 &   2.16 \\ \midrule 
 pglib\_opf\_case24\_ieee\_rts\_\_api & 13.01 & 6.89e+01  &  9.24 &   2.06 &  1.91e+01 &   8.27 \\ \midrule 
       pglib\_opf\_case30\_as\_\_api & 44.60 & 5.50e+01  & 12.60 &   1.41 &  1.74e+00 &   2.07 \\ \midrule 
      pglib\_opf\_case30\_fsr\_\_api &  2.75 & 2.14e+01  & 13.32 &   0.28 &  1.55e+00 &   1.01 \\ \midrule 
     pglib\_opf\_case30\_ieee\_\_api &  3.72 & 2.46e+01  &  5.80 &   0.00 &  7.22e-01 &   1.58 \\ \midrule 
     pglib\_opf\_case39\_epri\_\_api &  1.57 & 5.80e+01  &  8.25 &   0.18 &  8.18e+00 &   2.68 \\ \midrule 
     pglib\_opf\_case57\_ieee\_\_api &  0.07 & 2.15e+01  &  2.27 &   0.00 &  3.65e+00 &   0.94 \\ \midrule 
 pglib\_opf\_case73\_ieee\_rts\_\_api & 11.06 & 1.05e+02  &  7.37 &   2.91 &  2.14e+01 &   7.35 \\ \midrule 
   pglib\_opf\_case89\_pegase\_\_api &  8.13 & 4.02e+01  &  7.39 &   6.88 &  3.11e+01 &   6.26 \\ \midrule 
    pglib\_opf\_case118\_ieee\_\_api & 28.63 & 3.64e+02  & 16.46 &  11.14 &  3.30e+02 &  10.51 \\ \midrule 
pglib\_opf\_case162\_ieee\_dtc\_\_api &  5.44 & 1.03e+02  &  5.74 &   1.71 &  3.43e+02 &   4.72 \\ \midrule 
     pglib\_opf\_case179\_goc\_\_api &  7.17 & 2.38e+03  &  9.92 &   0.64 &  2.12e+03 &   8.96 \\ \midrule 
    pglib\_opf\_case200\_tamu\_\_api &  0.02 & 1.16e-01  &  0.68 &   0.00 &  4.10e-05 &   0.45 \\ \midrule 
   pglib\_opf\_case240\_pserc\_\_api &  0.79 & 4.19e+02  &  7.53 &   0.33 &  3.68e+02 &   7.16 \\ \midrule 
    pglib\_opf\_case300\_ieee\_\_api &  0.88 & 1.52e+02  &  3.81 &   0.03 &  3.79e+01 &   2.25 \\ \midrule 
    pglib\_opf\_case500\_tamu\_\_api &  0.07 & 1.98e+00  &  1.16 &   0.00 &  3.25e-03 &   0.48 \\ \midrule 
    pglib\_opf\_case588\_sdet\_\_api &  0.92 & 1.88e+03  &  5.32 &   0.53 &  2.39e+02 &   3.82 \\ \midrule 
 pglib\_opf\_case1354\_pegase\_\_api &  0.85 & 2.97e+04  &  5.55 &   0.38 &  3.76e+04 &   4.40 \\ \midrule 
    pglib\_opf\_case1888\_rte\_\_api &  0.28 & n.a.  &  6.44 &   0.17 &  n.a. &   5.84 \\  \midrule 
    pglib\_opf\_case1951\_rte\_\_api &  0.43 & n.a.  &  3.44 &   0.19 &  n.a. &   3.01 \\  \midrule 
   pglib\_opf\_case2316\_sdet\_\_api &  1.95 & 4.76e+03  &  7.34 &   0.61 &  1.16e+03 &   3.92 \\ \midrule 
    pglib\_opf\_case2383wp\_k\_\_api &  0.00 & 1.52e+03  &  4.97 &   0.00 &  3.43e+03 &   9.46 \\ \midrule 
    pglib\_opf\_case2736sp\_k\_\_api & 12.97 & 1.58e+03  &  5.98 &   2.60 &  1.57e+03 &   3.50 \\ \midrule 
   pglib\_opf\_case2737sop\_k\_\_api &  3.66 & 1.07e+03  &  9.87 &   2.83 &  1.10e+03 &   2.04 \\ \midrule 
   pglib\_opf\_case2746wop\_k\_\_api &  0.00 & 1.05e+04  &  4.52 &   0.00 &  1.58e+04 &   8.06 \\ \midrule 
    pglib\_opf\_case2746wp\_k\_\_api &  0.00 & 4.79e+03  &  4.46 &   0.00 &  7.82e+03 &   7.94 \\ \midrule 
    pglib\_opf\_case2848\_rte\_\_api &  0.21 & n.a.  &  3.71 &   0.05 &  n.a. &   2.87 \\  \midrule 
   pglib\_opf\_case2853\_sdet\_\_api &  2.30 & 4.26e+03  &  7.91 &   0.95 &  4.41e+03 &   4.59 \\ \midrule 
    pglib\_opf\_case2868\_rte\_\_api &  0.18 & n.a.  &  2.36 &   0.19 &  n.a. &   2.04 \\  \midrule 
 pglib\_opf\_case2869\_pegase\_\_api &  1.31 & 6.83e+02  &  2.84 &   0.92 &  5.98e+02 &   2.25 \\ \midrule 
    pglib\_opf\_case3012wp\_k\_\_api &  0.00 & 1.88e+03  &  5.56 &   0.00 &  3.79e+03 &   6.70 \\ \midrule 
    pglib\_opf\_case3120sp\_k\_\_api & 24.14 & 2.41e+03  &  8.13 &   9.59 &  5.16e+03 &   4.99 \\ 
\bottomrule
\multicolumn{7}{l}{n.a. -- The AC power flows in MATPOWER did not converge for these test cases.}\\
\end{tabular}
\end{table*}

\begin{table*}[]
    \centering
        \caption{Optimality gaps, distances to AC-feasibility and local optimality for PGLib OPF test cases from \cite{PGLIB}  (small phase angle difference cases (sad))}
    \label{Appendix_3}
\begin{tabular}{r c c c c c c}
\toprule
 & \multicolumn{3}{c}{QC relaxation} & \multicolumn{3}{c}{SDP relaxation} \\
\midrule
&  & distance to & distance to &  & distance to & distance to\\
 & opt. gap & AC-feas. & local opt.&  opt. gap  & AC-feas. & local opt. \\
Test cases & (\%) & (\%) & (\%) & (\%) &(\%) &(\%) \\
\midrule
pglib\_opf\_case14\_ieee\_\_sad &  7.16 & 3.82e+01  &  5.04 &   0.03 &  6.00e+00 &   4.56 \\ \midrule 
pglib\_opf\_case24\_ieee\_rts\_\_sad &  2.93 & 2.93e+01  &  5.71 &   2.52 &  1.48e+01 &  10.19 \\ \midrule 
      pglib\_opf\_case30\_as\_\_sad &  2.31 & 1.18e+01  &  7.65 &   0.16 &  9.39e+00 &   5.40 \\ \midrule 
     pglib\_opf\_case30\_fsr\_\_sad &  0.41 & 4.78e+00  &  9.03 &   0.02 &  1.70e+00 &   5.57 \\ \midrule 
    pglib\_opf\_case30\_ieee\_\_sad &  3.41 & 9.57e+00  &  2.72 &   0.00 &  4.80e-05 &   2.36 \\ \midrule 
    pglib\_opf\_case39\_epri\_\_sad &  0.19 & 9.63e+00  &  2.39 &   0.02 &  2.17e+01 &   5.54 \\ \midrule 
    pglib\_opf\_case57\_ieee\_\_sad &  0.83 & 1.17e+02  &  2.77 &   0.08 &  6.38e+01 &   4.96 \\ \midrule 
pglib\_opf\_case73\_ieee\_rts\_\_sad &  2.53 & 6.20e+01  &  4.72 &   1.47 &  6.69e+01 &   8.83 \\ \midrule 
  pglib\_opf\_case89\_pegase\_\_sad &  0.82 & 9.40e+01  &  5.44 &   0.05 &  3.67e+02 &   4.77 \\ \midrule 
   pglib\_opf\_case118\_ieee\_\_sad &  9.48 & 1.24e+02  &  8.12 &   3.70 &  7.00e+01 &   7.55 \\ \midrule 
sad/pglib\_opf\_case162\_ieee\_dtc\_\_sad &  8.02 & 3.10e+02  &  6.97 &   2.28 &  1.25e+03 &   7.66 \\ \midrule 
    pglib\_opf\_case179\_goc\_\_sad &  1.04 & 1.11e+03  &  8.03 &   0.95 &  2.56e+03 &  10.68 \\ \midrule 
   pglib\_opf\_case200\_tamu\_\_sad &  0.00 & 1.17e+00  &  0.51 &   0.00 &  2.42e-02 &   1.73 \\ \midrule 
  pglib\_opf\_case240\_pserc\_\_sad &  5.24 & 7.73e+02  & 10.32 &   4.17 &  7.33e+02 &  10.37 \\ \midrule 
   pglib\_opf\_case300\_ieee\_\_sad &  2.35 & 1.61e+02  &  2.64 &   0.11 &  1.87e+01 &   2.90 \\ \midrule 
   pglib\_opf\_case500\_tamu\_\_sad &  7.89 & 1.66e+01  &  1.33 &   7.59 &  1.72e+03 &   3.06 \\ \midrule 
   pglib\_opf\_case588\_sdet\_\_sad &  6.26 & 3.56e+03  &  9.49 &   5.75 &  6.57e+02 &   9.01 \\ \midrule 
pglib\_opf\_case1354\_pegase\_\_sad &  2.36 & 4.09e+03  &  3.63 &   0.59 &  1.76e+03 &   3.18 \\ \midrule 
   pglib\_opf\_case1888\_rte\_\_sad &  2.72 & n.a.  &  9.23 &   2.67 &  n.a. &  10.18 \\  \midrule 
   pglib\_opf\_case1951\_rte\_\_sad &  0.41 & n.a.  &  2.32 &   0.28 &  n.a. &   2.96 \\  \midrule 
  pglib\_opf\_case2316\_sdet\_\_sad &  2.17 & 6.93e+03  &  6.35 &   0.73 &  2.32e+03 &   5.03 \\ \midrule 
   pglib\_opf\_case2383wp\_k\_\_sad &  2.15 & 1.45e+04  &  5.26 &   0.56 &  9.59e+03 &   4.25 \\ \midrule 
   pglib\_opf\_case2736sp\_k\_\_sad &  1.53 & 1.74e+04  &  5.14 &   0.18 &  1.94e+04 &   2.90 \\ \midrule 
  pglib\_opf\_case2737sop\_k\_\_sad &  1.92 & 1.29e+04  &  4.28 &   0.51 &  2.22e+04 &   2.76 \\ \midrule 
  pglib\_opf\_case2746wop\_k\_\_sad &  2.00 & 4.03e+05  &  6.08 &   0.36 &  2.89e+04 &   2.65 \\ \midrule 
   pglib\_opf\_case2746wp\_k\_\_sad &  1.67 & 3.66e+04  &  6.00 &   0.19 &  1.80e+04 &   2.56 \\ \midrule 
   pglib\_opf\_case2848\_rte\_\_sad &  0.27 & n.a.  &  6.33 &   0.20 &  n.a. &   7.27 \\  \midrule 
  pglib\_opf\_case2853\_sdet\_\_sad &  2.39 & 4.95e+03  &  8.16 &   1.42 &  3.98e+03 &   6.86 \\ \midrule 
   pglib\_opf\_case2868\_rte\_\_sad &  0.54 & n.a.  &  2.34 &   0.51 &  n.a. &   4.41 \\  \midrule 
pglib\_opf\_case2869\_pegase\_\_sad &  1.41 & 9.49e+03  &  3.20 &   0.42 &  1.46e+03 &   2.90 \\ \midrule 
   pglib\_opf\_case3012wp\_k\_\_sad &  1.40 & 4.80e+03  &  4.39 &   0.45 &  1.22e+04 &   1.84 \\ \midrule 
   pglib\_opf\_case3120sp\_k\_\_sad &  1.41 & 7.83e+03  &  4.29 &   0.60 &  4.16e+03 &   2.98 \\ \midrule 
\bottomrule
\multicolumn{7}{l}{n.a. -- The AC power flows in MATPOWER did not converge for these test cases.} \\
\end{tabular}
\end{table*}

\end{document}